\documentclass{amsart}
\usepackage{verbatim}
\usepackage{amssymb}
\usepackage{amsbsy}
\usepackage{amscd}
\usepackage{amsmath}
\usepackage{amsthm}
\usepackage[mathscr]{eucal}

\newenvironment{NB}{
{\bf NB}. \footnotesize
}{}
\renewenvironment{NB}{
\comment
  }{\endcomment}

\theoremstyle{plain}
 \newtheorem{thm}{Theorem}[section]
 \newtheorem{lem}[thm]{Lemma}
 \newtheorem{prop}[thm]{Proposition}
 \newtheorem{cor}[thm]{Corollary}
 
\theoremstyle{definition}
 \newtheorem{defn}{Definition}[section]
\theoremstyle{remark}
 \newtheorem{rem}{Remark}[section]
 \newtheorem{ex}{Example}[section]
 \newtheorem{claim}{Claim}[section]
\def\Bbb{\mathbb}
\def\frak{\mathfrak}
\def\cal{\mathcal}

\newcommand{ \Supp}{\operatorname{Supp}}
\newcommand{\Ext}{\operatorname{Ext}}
\newcommand{\Hom}{\operatorname{Hom}}
\newcommand{\End}{\operatorname{End}}

\newcommand{\im}{\operatorname{im}}

\newcommand{\rk}{\operatorname{rk}}

\newcommand{\NS}{\operatorname{NS}}
\newcommand{\coker}{\operatorname{coker}}
\newcommand{\Pic}{\operatorname{Pic}}
\newcommand{\ch}{\operatorname{ch}}
\newcommand{\td}{\operatorname{td}}

\newcommand{\Hilb}{\operatorname{Hilb}}

\newcommand{\Coh}{\operatorname{Coh}}

\newcommand{\Div}{\operatorname{Div}}

\newcommand{\id}{\operatorname{id}}

\newcommand{\tr}{\operatorname{tr}}
\newcommand{\ad}{\operatorname{ad}}
\newcommand{\Fun}{\operatorname{Fun}}

\font\b=cmr10 scaled \magstep5
\def\bigzerou{\smash{\lower1.7ex\hbox{\b 0}}}

\numberwithin{equation}{section}
\begin{document}

\title
{An action of a Lie algebra on 
the homology groups of moduli spaces of stable sheaves}
\author{K\={o}ta Yoshioka\
}
 \address{Department of Mathematics, Faculty of Science,
Kobe University,
Kobe, 657, Japan}
\email{yoshioka@@math.kobe-u.ac.jp}

\begin{abstract}
We construct an action of a Lie algebra on the homology groups
of moduli spaces of stable sheaves on $K3$ surfaces 
under some technical conditions.
This is a generalization of Nakajima's construction of $\frak{sl}_2$-action
on the homology groups \cite{N:2002}.
In particular, for an $A,D,E$-configulation of $(-2)$-curves,
we shall give a collection of moduli spaces
such that the associated Lie algebra acts on their homology groups.
\end{abstract} 
\maketitle

\section{Introduction}

Let $X$ be a smooth projective surface defined over ${\Bbb C}$
and $H$ an ample divisor on $X$.
Assume that $X$ is a $K3$ surface.
Let $M_H(v)$ be the moduli space of $H$-stable sheaves $E$  
with the Mukai vector $v(E)=v$ (cf. \eqref{eq:mukai-vector}). 
In \cite{Y:5}, we studied 
a special kind of Fourier-Mukai transform called
$(-2)$-reflection. 
For this purpose, we introduced the Brill-Noether locus on the moduli
space and studied its properties. 
Similar results are obtained by Markman \cite{Mar}.
We fix a vector bundle $G$ on $X$.
A stable sheaf $E_0$ is called exceptional, if
$\Ext^1(E_0,E_0)=0$.
Then $v(E_0)$ is a $(-2)$-vector, that is,
$\langle v(E_0)^2 \rangle=-2$.
We assume that the twisted degree 
$\deg_G(E_0):=\deg(G^{\vee} \otimes E_0)=0$.
Let $v \in H^*(X,{\Bbb Z})$ be a Mukai vector such that
\begin{equation}\label{eq:bn-condition}
 \deg_G(E)=\min\{\deg_G(E')>0| E' \in K(X) \}
\end{equation}
for $E \in M_H(v)$.
Let 
\begin{equation}
M_H(v)_{E_0,n}:=
\{E \in M_H(v)|\dim \Hom(E_0,E)=n \}
\end{equation}
be {\it the Brill-Noether locus} with respect to $E_0$.
Under the condition \eqref{eq:bn-condition}, 
we showed that $M_H(v)_{E_0,n}$ is a Grassmannian bundle over a
smooth manifold such that the relative cotangent bundle is
isomorphic to the normal bundle $N_{M_H(v)_{E_0,n}/M_H(v)}$.
Similar Grassmannian structure appears in Nakajima's quiver varieties
\cite{Na:1998}. By usnig this structure, he constructed a Lie algebra
action on the (Borel-Moore) homology groups of quiver varieties.
Based on our description of the Brill-Noether locus,
recently Nakajima \cite{N:2002} constructed an ${\frak {sl}}_2$-action
on the homology groups of moduli spaces
$\bigoplus_v H_*(M_H(v),{\Bbb C})$,
where $v$ runs a suitable set of Mukai vectors satisfing
minimality condition \eqref{eq:bn-condition}.
  
In this note, under the same condition,
we shall generalize Nakajima's result.
Thus we shall construct a Lie algebra action on
the homology groups of moduli spaces of stable sheaves
(Theorem \ref{thm:action}):
For a collection of exceptional sheaves $E_i$, $i=1,2,\dots,s$ which
satisfy some technical conditions,
we shall construct operators $h_i,e_i,f_i$, $i=1,2,\dots,s$
and show that they satisfy the commutation relations for
Chevalley generators.
In particular, we show that $[e_i,f_i]=h_i$ and 
$[e_i,f_j]=0$, $i \ne j$.
Since the first relation is proved by Nakajima,
we only need to show the second one.
For this purpose, we introduce the notion of 
universal extension (resp. division) 
with respect to $E_i$, $i=1,2,\dots, s$
(see, sect. \ref{subsect:univ}).
This is our main idea and the other arguments are included in
Nakajima's papers. 
Since the action is defined by algebraic correspondences,
we also have an action on the rational Chow groups. 
In section \ref{sect:example},
we give some examples of actions.

 
Replacing $E_0$ by a purely 1-dimensional exceptional sheaf  
and the minimality condition by 
$\chi(E)=1$, our construction also works for moduli spaces
of purely 1-dimensional stable sheaves.
In particular, we shall construct an action of the affine Lie algebra
associated to a singular fiber of an elliptic surface.  
On an elliptic surface, purely 1-dimensional sheaves are related to 
torsion free sheaves of relative degree 0 via 
the relative Fourier-Mukai transform.
Moreover purely 1-dimensional sheaves are related to
the enumerative geometry of curves on $X$
(cf. \cite{Y-Z:1}). 
Thus the moduli spaces of purely 1-dimensional stable sheaves
are important objects to study.
For a rational elliptic surface $X$, it is observed in \cite{MNWV}
that
the Euler characteristics of the moduli spaces are $W(E_8^{(1)})$-invariant,
where $W(E_8^{(1)})$ is the Weyl group
associated to the $E_8^{(1)}$-lattice $K_X^{\perp} \subset H^2(X,{\Bbb Z})$.
An explanation is given in terms of the monodromy action,
that is, we use the invariance 
of the homology groups of the moduli spaces under the 
deformation of $X$.
Our construction of the Lie algebra gives another explanation 
of this invariance. These are treated in section \ref{sect:1-dim}.
In section \ref{sect:equiv}, we give a remark 
on the case of $G$-equivariant sheaves.

\section{Moduli of stable sheaves of minimal degree}\label{sect:moduli}

{\it Notation.}

Let $X$ be a smooth projective surface.
Let $\Coh(X)$ be the category of coherent sheaves on $X$
and $K(X)$ the Grothendieck group of $X$.
In this paper, we use the Borel-Moore homology groups.
For an algebraic set $M$,
$H_*(M,{\Bbb C})$ denotes the Borel-Moore homology group
of $M$.
If $M$ is compact, then $H_*(M,{\Bbb C})$ coincides with the
usual singular homology group of $M$.

Let 
${\bf D}(X):={\bf D}^b(\Coh(X))$ be the bounded derived categories of
$\Coh(X)$.
For complexes ${\Bbb E}, {\Bbb F} \in {\bf D}(X)$,
we set $\Ext^i({\Bbb E},{\Bbb F}):=
\Hom_{{\bf D}(X)}({\Bbb E},{\Bbb F}[i])$.
We usually denote $\Ext^0({\Bbb E},{\Bbb F})$ by
$\Hom({\Bbb E},{\Bbb F})$.
For a morphism $\phi:{\Bbb E} \to {\Bbb F}$,
$[{\Bbb E},{\Bbb F}]$ denotes the mapping cone of a representative of
$\phi$.
If $H^i([{\Bbb E} \to {\Bbb F}])=0$ for all $i$, then we write
${\Bbb E} \cong {\Bbb F}$.
We usually denotes $\Ext^i([{\Bbb E}_1 \to {\Bbb E}_2],{\Bbb F})$
(resp. $\Ext^i({\Bbb F},[{\Bbb E}_1 \to {\Bbb E}_2])$)
by $\Ext^i({\Bbb E}_1 \to {\Bbb E}_2,{\Bbb F})$
(resp. $\Ext^i({\Bbb F},{\Bbb E}_1 \to {\Bbb E}_2)$).

Let $H$ be an ample divisor on $X$ and $G$ an element of $K(X)$ 
with $\rk G>0$.
For a coherent sheaf $E$ on $X$, we set
$\deg_G(E):=\deg(G^{\vee} \otimes E)$ and
$\chi_G(E):=\chi(G^{\vee} \otimes E)$.

\subsection{Technical lemmas}\label{subsect:tech}

In this subsection, we introduce some technical conditions
\eqref{eq:min-cond}, \eqref{eq:min-cond2}, \eqref{eq:min-cond2'}
and under these conditions
we give some technical lemmas.
These will play important roles for
our construction of the action. 

\begin{defn}
A purely 1-dimensional sheaf $E$ is $\mu$-stable,
if the scheme-theoretic support $\Div(E)$ of $E$ is reduced and irreducible.
\end{defn}
We fix an ample divisor $H$ on $X$.
Let $G$ be an element of $K(X)$ with $\rk G>0$.
In this note, we treat $\mu$-semi-stable sheaves $E$ with
\begin{equation}\label{eq:min-cond}
\deg_G(E)=\min\{\deg_G(E')>0 |\;E' \in K(X) \}.
\end{equation}
This is a fairly strong condition for $E$, but
such $E$ behave very well.
\begin{lem}\label{lem:key1}
Let $G$ be an element of $K(X)$ with $\rk G>0$ and $E_i$,
$i=1,2,\dots,s$ be 
$\mu$-stable vector bundles with $\deg_G(E_i)=0$.
Let $E$ be a $\mu$-semi-stable sheaf satisfying \eqref{eq:min-cond}.

\begin{enumerate}
\item[(1)]
Then $E$ is $\mu$-stable.
\item[(2)]
Every non-trivial extension
\begin{equation}
0 \to E_1 \to F \to E \to 0
\end{equation}
defines a $\mu$-stable sheaf.
\item[(3)]
Let $V_i$ be subspaces of $\Hom(E_i,E)$, $i=1,2,\dots,s$.
Then  
$\phi:\bigoplus_{i=1}^s V_i \otimes E_i \to E$ 
is injective or surjective in codimension 1.
Moreover,
\begin{enumerate}
\item[(3-1)] 
if $\phi:\bigoplus_{i=1}^s V_i \otimes E_i \to E$ is injective, then the
cokernel is $\mu$-stable,
\item[(3-2)]
if $\phi:\bigoplus_{i=1}^s V_i \otimes E_i \to E$ 
is surjective in codimension 1, then
$\ker \phi$ is $\mu$-stable.
In particular $D(E):={\cal E}xt^1(\bigoplus_{i=1}^s V_i \otimes E_i \to E,{\cal O}_X)$
is $\mu$-stable.
\end{enumerate}
\end{enumerate}
\end{lem}
Since $\deg_G(E)/\rk(E)=\rk (G)(\deg(E)/\rk E-\deg(G)/\rk G)$,
the $\mu$-stability can be defined by using the $G$-twisted slope
$\deg_G(E)/\rk(E)$.
By using the following lemmas, the proof of \cite[Lem. 2.1]{Y:5}
implies our lemma.
So we only give a proof of (1), (3).
We first note the following easy lemmas.

\begin{lem}\label{lem:integral}
A purely 1-dimensional sheaf $E$ with
\eqref{eq:min-cond} is $\mu$-stable.
\end{lem}

\begin{lem}\label{lem:ineq}
Let $r,d,x$ be positive integers.
Let $y$ be an integer such that $y \in d {\Bbb Z}$.  
If $0<y/x<d/r$, then $y \geq d$ and $x>r$.
\end{lem}

{\it Proof of Lemma \ref{lem:key1} (1), (3).}
Let $E'$ be a subsheaf of $E$ with $\deg_G(E)/\rk E=\deg_G(E')/\rk E'$.
Then $1 \geq \deg_G(E)/\deg_G(E')=\rk E/\rk E' \geq 1$.
Hence $\rk E'=\rk E$ and $\deg_G(E')=\deg_G(E)$, which implies that
$E$ is $\mu$-stable.
Thus (1) holds. We shall prove (3).
We first assume that $\rk E>0$.
By the $\mu$-stability of $E$, we have
\begin{equation}
0 \leq \frac{\deg_G(\im \phi)}{\rk (\im \phi)} \leq 
\frac{\deg_G(E)}{\rk E}.
\end{equation}
By Lemma \ref{lem:ineq},
(i) $\deg_G(\im \phi)=0$ or (ii) 
$\deg_G(\im \phi)/\rk (\im \phi)=\deg_G(E)/\rk E$.
In the first case, $\deg_G(\ker \phi)=0$.
Assume that $\ker \phi \ne 0$.
Let $F$ be a $\mu$-stable locally free subsheaf of $\ker \phi$ with
$\deg_G(F)=0$.
Then there is a non-zero homomorphism $F \to E_1$, which is isomorphic.
Hence $\Hom(E_1,\ker \phi) \ne 0$, which is a contradiction.
Therefore $\ker \phi=0$.
We shall show that $E':=\coker \phi$ is $\mu$-stable.
We note that $E'$ does not have a 0-dimensional subsheaf and
$\deg_G(E')=\deg_G(E)$.
We first assume that $\rk E'>0$.
If $E'$ is not $\mu$-stable, then (1) implies that
$E'$ is not $\mu$-semi-stable.
Then there is a quotient $E' \to F$ with $\deg_G(F)/\rk F<\deg_G(E')/\rk E'$.
By Lemma \ref{lem:ineq}, $\deg_G(F) \leq 0$, which implies that
$\deg_G(F)/\rk F<\deg_G(E)/\rk E$.
This is a contradiction.
Therefore $E'$ is $\mu$-stable.
If $\rk E'=0$, then $E'$ is of pure dimension 1.
Then Lemma \ref{lem:integral} implies that $E'$ is $\mu$-stable.

We next treat the second case:
$\deg_G(\im \phi)/\rk (\im \phi)=\deg_G(E)/\rk E$.
In this case, $\phi$ is surjective in codimension 1.
We shall show that $\ker \phi$ is $\mu$-stable.
Assume that there is a locally free subsheaf $F$ of $\ker \phi$
with $\deg_G(F)/\rk F>\deg_G(\ker \phi)/\rk (\ker\phi)=
-\deg_G(E)/\rk (\ker \phi)$.
Then we get that $\deg_G(F) \leq 0$.
If $\deg_G(F)=0$, then $\Hom(F,E_i) \ne 0$ for an $i$.
Since $E_i$ and $F$ are $\mu$-stable sheaves with the same slope,
non-trivial homomorphism
$F \to E_i$ is isomorphic in codimension 1.
Since $F$ is locally free, we conclude that $F \cong E_i$.
Then $\Hom(E_i,F) \ne 0$, which is a contradiction.
Hence $\deg_G(F)<0$, which means that
$0<-\deg_G(F)/\rk F<\deg_G(E)/\rk(\ker \phi)$,
Then
Lemma \ref{lem:ineq} implies that
$-\deg_G(F)\geq \deg_G(E)$ and $\rk F>\rk (\ker \phi)$,
which is a contradiction.
Therefore $\ker \phi$ is $\mu$-stable.

If $\rk E=0$, then since $E$ is $\mu$-stable, we get
$\phi=0$ or $\phi$ is surjective in codimension 1.
Then by the same arguments as above, we see that 
$\ker \phi$ is $\mu$-stable.   
\qed

Besides the condition for $\mu$-semi-stable sheaves \eqref{eq:min-cond},
we also introduce similar conditions and lemmas for Gieseker
(twisted) semi-stabilities. 
\begin{defn}
Let $G$ be an element of $K(X)$ with $\rk G>0$. 
A torsion free sheaf $E$ is $G$-twisted stable, if 
\begin{equation}
\frac{\chi_G(F(nH))}{\rk F}<\frac{\chi_G(E(nH))}{\rk E}, n \gg 0
\end{equation}
for all proper subsheaf $F (\ne 0)$ of $E$.
\end{defn} 
As in the proof of Lemma \ref{lem:key1},
we also have the following assertions.

\begin{lem}\label{lem:key2}
Let $G$ be an element of $K(X)$ with $\rk G>0$ and 
$E_i$, $i=1,2,\dots,s$, be
$G$-twisted stable sheaves with $\deg_G(E_i)=\chi_G(E_i)=0$.
Let $E$ be a $G$-twisted stable torsion free sheaf with 
$\deg_G(E)=0$ and
\begin{equation}\label{eq:min-cond2}
\chi_G(E)=\min\{\chi_G(E')>0|\;E' \in \Coh(X), \deg_G(E')=0 \}
\end{equation}
or $E={\Bbb C}_P$, $P \in X$ with \eqref{eq:min-cond2}.
\begin{enumerate}
\item[(1)]
Then every non-trivial extension
\begin{equation}
0 \to E_1 \to F \to E \to 0
\end{equation}
defines a $G$-twisted stable sheaf.
\item[(2)]
Let $V_i$ be a subspace of $\Hom(E_i,E)$.
Then  
$\phi:\bigoplus_{i=1}^s V_i \otimes E_i \to E$ 
is injective or surjective. Moreover,
\begin{enumerate}
\item[(2-1)] 
if $\phi:\bigoplus_{i=1}^s V_i \otimes E_i \to E$ is injective, then the
cokernel is a $G$-twisted stable torsion free sheaf or
${\Bbb C}_P, P \in X$,
\item[(2-2)]
if $\phi:\bigoplus_{i=1}^s V_i \otimes E_i \to E$ is surjective, then
$\ker \phi$ is $G$-twisted stable.
\end{enumerate}
\end{enumerate}
\end{lem}

\begin{lem}\label{lem:key2'}
Let $G$ be an element of $K(X)$ with $\rk G>0$ and $E_i$,
$i=1,2,\dots,s$, be $G$-twisted stable sheaf with 
$\deg_G(E_i)=\chi_G(E_i)=0$.
Let $E$ be a $G$-twisted stable torsion free sheaf with 
$\deg_G(E)=0$ and
\begin{equation}\label{eq:min-cond2'}
\chi_G(E)=\max\{\chi_G(E')<0|\;E' \in \Coh(X), \deg_G(E')=0 \}.
\end{equation}
\begin{enumerate}
\item[(1)]
Then every non-trivial extension
\begin{equation}
0 \to E \to F \to E_1 \to 0
\end{equation}
defines a $G$-twisted stable sheaf.
\item[(2)]
Let $V_i$ be a subspace of $\Hom(E,E_i)$.
Then  
$\phi:E \to \bigoplus_{i=1}^s V_i^{\vee} \otimes E_i $ 
is injective or surjective.
Moreover,
\begin{enumerate}
\item[(2-1)] 
if $\phi:E \to \bigoplus_{i=1}^s V_i^{\vee} \otimes E_i$ 
is injective, then the
cokernel is a $G$-twisted stable torsion free sheaf or
${\Bbb C}_P$, $P \in X$
\item[(2-2)]
if $\phi:E \to \bigoplus_{i=1}^s V_i^{\vee} \otimes E_i$ is surjective, 
then
$\ker \phi$ is $G$-twisted stable.
\end{enumerate}
\end{enumerate}
\end{lem}

\subsection{Basic properties of stable sheaves of minimal degree}

Assume that $K_X$ is numerically trivial.
We define a bilinear form $\langle \quad,\quad \rangle$
on $H^*(X,{\Bbb Q}):=\bigoplus_{i=0}^2
H^{2i}(X,{\Bbb Q})$ by
\begin{equation}\label{eq:mukai-lattice}
\langle x,y \rangle:=\int_X x_1 \wedge y_1-x_0 \wedge y_2-x_2 \wedge y_0
\end{equation}
where $x_i \in H^{2i}(X,{\Bbb Q})$ (resp. $y_i \in H^{2i}(X,{\Bbb Q})$)
is the $2i$-th component of $x$ (resp. $y$).

For an object ${\Bbb E} \in {\bf D}(X)$,
we define the Mukai vector of ${\Bbb E}$ by
\begin{equation}\label{eq:mukai-vector}
\begin{split}
v({\Bbb E})&=\sum_i (-1)^i v(H^i({\Bbb E}))\\
&=\sum_i (-1)^i \ch(H^i({\Bbb E}))\sqrt{\td_X} \in H^*(X,{\Bbb Q}),
\end{split}
\end{equation}
where $\td_X$ is the todd class of $X$.
We have a map $v:{\bf D}(X) \to H^*(X,{\Bbb Q})$.
We call an element of $v({\bf D}(X))$ a Mukai vector.
For ${\Bbb E}, {\Bbb F} \in {\bf D}(X)$,
we define the Riemann-Roch number by
\begin{equation}
\chi({\Bbb E},{\Bbb F}):=\sum_i (-1)^i\dim \Ext^i({\Bbb E},{\Bbb F}).
\end{equation} 
Then the Riemann-Roch theorem says the following.
\begin{prop}
\begin{equation}
\chi({\Bbb E},{\Bbb F})=-\langle v({\Bbb E}),v({\Bbb F}) \rangle.
\end{equation}
\end{prop}
By a similar way, we also define the rank $\rk {\Bbb E}$
and other invariants.  
We fix an element $G \in K(X)$ with $\rk G>0$.
For an object ${\Bbb E} \in {\bf D}(X)$
such that
$\deg_G({\Bbb E})$ satisfies \eqref{eq:min-cond},
we define a stability condition.
\begin{defn}\label{defn:stability}
Let ${\Bbb E} \in {\bf D}(X)$ be an object such that
$\deg_G({\Bbb E})$ satisfies \eqref{eq:min-cond}.
Then ${\Bbb E}$ is stable, if 
\begin{equation}\label{eq:length-cond}
H^{i}({\Bbb E} \overset{\Bbb L}{\otimes}
 {\Bbb C}_P)=0, i \ne -1,0
\end{equation}
for all $P \in X$ and
one of the following conditions holds:
\begin{enumerate}
\item
$H^{i}({\Bbb E})=0$, $i \ne 0$ and $H^0({\Bbb E})$ is a stable sheaf.
\item
$H^i({\Bbb E})=0$, $i \ne -1,0$, 
$H^{-1}({\Bbb E})^{\vee}$ is a stable sheaf and
$H^0({\Bbb E})$ is a 0-dimensional sheaf.
\end{enumerate}
\end{defn}

\begin{rem}\label{rem:dual}
\begin{enumerate}
\item
The condition \eqref{eq:length-cond} implies that
there is a complex $C_{-1} \to C_0$ of locally free sheaves
which is quasi-isomorphic to ${\Bbb E}$. 
\item
If $\rk {\Bbb E}<0$, then $H^i(D({\Bbb E}))=0$,
$i \ne 1$ and $H^1(D({\Bbb E}))$ is a stable sheaf,
where $D({\Bbb E}):={\bf R}{\cal H}om({\Bbb E},{\cal O}_X)$
is the dual of ${\Bbb E}$.
Since we want to treat two cases simultaniously,
we use ${\Bbb E}$ instead of using $D({\Bbb E})$. 
\end{enumerate}
\end{rem}

\begin{defn}
For a Mukai vector $v \in H^*(X,{\Bbb Q})$ with the property
\eqref{eq:min-cond},
let $M_H(v)$ be the moduli space of 
(quasi-isomorphism classes of) stable complexes ${\Bbb E}$
with $v({\Bbb E})=v$. 
\end{defn}
If $\rk v<0$, then by Remark \ref{rem:dual},
$M_H(v)$ has a scheme structure.
The Zariski tangent space of $M_H(v)$ at ${\Bbb E}$ is
$\Ext^1({\Bbb E},{\Bbb E})$
and the obstruction for the infinitesimal liftings belongs to
the kernel of the trace map
\begin{equation}
\tr:\Ext^2({\Bbb E},{\Bbb E}) \to H^2(X,{\cal O}_X).
\end{equation}
In this paper, we assume that the trace map
\begin{equation}\label{eq:smooth-condition}
\tr:\Ext^2({\Bbb E},{\Bbb E}) \to H^2(X,{\cal O}_X)
\end{equation}
is isomorphic.

By Lemma \ref{lem:key1} and the condition \eqref{eq:smooth-condition}, 
we get the following assertions.
\begin{lem}\label{lem:non-empty}
Assume that $v \in H^*(X,{\Bbb Q})$ satisfies \eqref{eq:min-cond}.
\begin{enumerate}
\item
If $M_H(v) \ne \emptyset$, then $\dim M_H(v)=\langle v^2 \rangle+1+p_g$.
In particular,
if there is a stable complex ${\Bbb E}$ with 
$v({\Bbb E})=v$, then
$\langle v({\Bbb E})^2 \rangle \geq -(p_g+1)$.
\item
Assume that $X$ is a K3 surface.
Then there is a stable complex ${\Bbb E}$
with $v({\Bbb E})=v$ if and only if
$\langle v^2 \rangle \geq -2$.
\end{enumerate}
\end{lem}
Fot the proof of (ii), we also use \cite[Thm. 0.2]{Y:5}.

Let $S:=\{E_1,E_2,\dots,E_n \}$ be a finite set of $\mu$-stable vector bundles 
such that $\deg_G(E_i)=0$, $1 \leq i \leq n$.
We assume that 
\begin{equation}\label{eq:(-2)-condition1}
E_i \otimes K_X \cong E_i,\; E_i \in S.
\end{equation}
Let ${\cal S}$ be a subcategory of $\Coh(X)$ 
consisting of semi-stable sheaves $F$ whose Jordan-H\"{o}lder grading
is $\bigoplus_i E_i^{\oplus n_i}$.

\begin{lem}\label{lem:hom[-1]}
$\Hom({\Bbb E},F)=0$ and $\Hom(F[1],{\Bbb E})=0$ for $F \in {\cal S}$.
\end{lem}

\begin{proof}
We use the spectral sequence
\begin{equation}
\begin{split}
E_2^{p,q}&=\bigoplus_{q'+q''=q}\Ext^p(H^{-q'}(*),H^{q''}(**))\\
& \Longrightarrow 
E_{\infty}^{p+q}=\Ext^{p+q}(*,**).
\end{split}
\end{equation}
Since $H^i({\Bbb E})=0$, $i \ne -1,0$,
$\Hom({\Bbb E},F)=\Hom(H^0({\Bbb E}),F)$.
If $\rk {\Bbb E} \geq 0$, then $H^0({\Bbb E})$ is a stable sheaf
of positive $G$-twisted degree.
Hence $\Hom(H^i({\Bbb E}),F)=0$.
If $\rk {\Bbb E} < 0$, then $H^0({\Bbb E})$ is a 0-dimension sheaf.
Hence $\Hom(H^0({\Bbb E}),F)=0$.
Therefore the first claim holds.
Since $\Hom(F[1],{\Bbb E})=\Hom(F,H^{-1}({\Bbb E}))$, we also get the second
claim.
\end{proof}

\subsection{A universal division and a universal extension}
\label{subsect:univ}

\begin{defn}
An exact triangle 
\begin{equation}
F \to {\Bbb E} \to \widetilde{\Bbb E} \to F[1]
\end{equation}
is a universal division of ${\Bbb E}$ with respect to
$\{E_1,E_2,\dots,E_n \}$, if
$F \in {\cal S}$
and $\widetilde{\Bbb E}$ is a stable complex such that
$\Hom(E_i,\widetilde{\Bbb E})=0$, $1 \leq i \leq n$. 
\end{defn}

For an exact triangle
\begin{equation}
F' \to {\Bbb E} \to {\Bbb E}' \to F'[1],
\end{equation}
we have an exact sequence
\begin{equation}
\Hom(F'[1],\widetilde{\Bbb E}) \to
\Hom({\Bbb E}', \widetilde{\Bbb E}) \to
\Hom({\Bbb E}, \widetilde{\Bbb E}) \to 
\Hom(F',\widetilde{\Bbb E}).
\end{equation}
By our assumption and Lemma \ref{lem:hom[-1]},
$\Hom({\Bbb E}', \widetilde{\Bbb E}) \to
\Hom({\Bbb E}, \widetilde{\Bbb E})$ is an isomorphism.
Hence we have a unique morphism 
${\Bbb E}' \to \widetilde{\Bbb E}$ in ${\bf D}(X)$
which induces a commutative diagram of exact triangles
(in ${\bf D}(X)$):
\begin{equation}
\begin{CD}
F' @>>> {\Bbb E} @>>> {\Bbb E}' @>>> F'[1]\\
@VVV @| @VVV @VVV \\
F @>>> {\Bbb E} @>>> \widetilde{\Bbb E} @>>> F[1].\\
\end{CD}
\end{equation}
In particular, a universal division of ${\Bbb E}$
is unique (up to isomorphism class in ${\bf D}(X)$). 
Since $\Hom(\widetilde{\Bbb E},\widetilde{\Bbb E}) \cong {\Bbb C}$ and
$\Hom(F,\widetilde{\Bbb E})=0$, we get
\begin{equation}\label{eq:univ-div[hom]}
\Hom({\Bbb E},\widetilde{\Bbb E}) \cong {\Bbb C}.
\end{equation}
Since $E_i \otimes K_X \cong E_i$, we see that
$\Hom(F \otimes K_X^{\vee},\widetilde{\Bbb E})=
\Hom(F \otimes K_X^{\vee}[1],\widetilde{\Bbb E})=0$.
Hence we also get that
\begin{equation}\label{eq:univ-div[hom^2]}
\Hom({\Bbb E},\widetilde{\Bbb E} \otimes K_X) \cong
\Hom(\widetilde{\Bbb E},\widetilde{\Bbb E} \otimes K_X)
\cong H^0(X,K_X).
\end{equation}

\begin{defn}
An exact triangle 
\begin{equation}
F \to \widehat{\Bbb E} \to {\Bbb E} \to F[1]
\end{equation}
is a universal extension of ${\Bbb E}$ with respect to
$\{E_1,E_2,\dots,E_n \}$, if
$F \in {\cal S}$ and $\widehat{\Bbb E}$ is a stable complex such that
$\Ext^1(\widehat{\Bbb E},E_i)=0$, $1 \leq i \leq n$. 
\end{defn}

For an exact triangle
\begin{equation}
F' \to {\Bbb E}' \to {\Bbb E} \to F'[1],
\end{equation}
we have an exact sequence
\begin{equation}
\Hom(\widehat{\Bbb E},F') \to
\Hom(\widehat{\Bbb E}, {\Bbb E}') \to
\Hom(\widehat{\Bbb E}, {\Bbb E}) \to 
\Ext^1(\widehat{\Bbb E},F').
\end{equation}
By our assumption and Lemma \ref{lem:hom[-1]},
$\Hom(\widehat{\Bbb E}, {\Bbb E}') \to
\Hom(\widehat{\Bbb E},{\Bbb E})$ is an isomorphism.
Hence we have a unique morphism $\widehat{\Bbb E} \to {\Bbb E}'$
which induces a commutative diagram of exact triangles
\begin{equation}
\begin{CD}
F' @>>> {\Bbb E}' @>>> {\Bbb E} @>>> F'[1]\\
@AAA @| @AAA @AAA \\
F @>>> \widehat{\Bbb E} @>>> {\Bbb E} @>>> F[1].\\
\end{CD}
\end{equation}
In particular, a universal extension of ${\Bbb E}$
is unique. 
For a universal extension, we also see that 
\begin{equation}\label{eq:uni-ext[hom]}
 \begin{split}
  \Hom(\widehat{\Bbb E},{\Bbb E})& \cong{\Bbb C},\\
  \Hom(\widehat{\Bbb E},{\Bbb E} \otimes K_X) & \cong
  \Hom(\widehat{\Bbb E},\widehat{\Bbb E} \otimes K_X) \cong H^0(X,K_X).
 \end{split}
\end{equation}

\subsubsection{Condition for the existence}

\begin{lem}\label{lem:univ-exist}
\begin{enumerate}
\item
If $S$ defines a negative definite lattice, then
a universal extension and a universal division exist
for ${\Bbb E}$.
\item
Assume that $S$ defines a negative semi-definite lattice of affine type.
Let $\delta:=\sum_i a_i v(E_i)$ satisfy $\langle \delta,v(E_i) \rangle=0$
for all $E_i \in S$.
If $\langle v({\Bbb E}),\delta \rangle \ne 0$, then
a universal extension or a universal division exist
for ${\Bbb E}$.
\end{enumerate}
\end{lem}

\begin{proof}
\begin{claim}\label{claim:1}
For a non-zero morphism $\psi:E_{n_1} \to {\Bbb E}$,
${\Bbb E}^{(1)}:=[E_{n_1} \to {\Bbb E}]$ is also stable.
\end{claim}
Proof of Claim \ref{claim:1}:
For a non-zero morphism $\psi:E_{n_1} \to {\Bbb E}$,
we have an exact sequence
\begin{equation}
 \Hom(E_{n_1},{\Bbb E}[-1]) \to
 \Hom(E_{n_1},{\Bbb E}^{(1)}[-1]) \to {\Bbb C} \overset{\psi}{\to}
 \Hom(E_{n_1},{\Bbb E}).
\end{equation}
By Lemma \ref{lem:hom[-1]},
$\Hom(E_{n_1},{\Bbb E}[-1])=0$, and hence we get
$\Hom(E_{n_1},{\Bbb E}^{(1)}[-1])=0$.
We note that ${\Bbb E}^{(1)}$ satisfies \eqref{eq:length-cond}
and we have the following exact sequence
\begin{equation}
 \begin{CD}
  0 @>>> H^{-1}({\Bbb E}) @>>> H^{-1}({\Bbb E}^{(1)}) @>>> @.\\
  E_{n_1} @>>> H^0({\Bbb E}) @>>> H^0({\Bbb E}^{(1)}) @>>> 0.    
 \end{CD}
\end{equation}
Then we get $0=\Hom(E_{n_1},{\Bbb E}^{(1)}[-1])\cong
\Hom(E_{n_1},H^{-1}({\Bbb E}^{(1)}))$.
If $H^{-1}({\Bbb E})=0$, then
$E_{n_1} \to H^0({\Bbb E})$ is a non-zero homomorphism.
\begin{NB}
If it is zero, then 
$\Hom(E_{n_1},H^{-1}({\Bbb E}^{(1)})) \ne 0$.
\end{NB}
By Lemma \ref{lem:key1},
${\Bbb E}^{(1)}$ is stable. 
Assume that $H^{-1}({\Bbb E}) \ne 0$.
Since $H^{-1}({\Bbb E}^{(1)})$ is locally free and
$\Hom(E_{n_1},H^{-1}({\Bbb E}^{(1)}))=0$,
the extension
\begin{equation}
 0 \to H^{-1}({\Bbb E}) \to  H^{-1}({\Bbb E}^{(1)}) \to
 E_{n_1}' \to 0
\end{equation}
does not split, where $E_{n_1}'$ is a subsheaf of
$E_{n_1}$ with $(E_{n_1}')^{\vee \vee}=E_{n_1}$.
By Lemma \ref{lem:key1},
${\Bbb E}^{(1)}$ is a stable complex. 
Thus the claim holds.

If there is a non-zero morphism $E_{n_2} \to {\Bbb E}^{(1)}$,
then we set ${\Bbb E}^{(2)}:=[E_{n_2} \to {\Bbb E}^{(1)}]$.
Then we have an exact triangle
\begin{equation}
F^2 \to {\Bbb E} \to {\Bbb E}^{(2)} \to F^2[1],
\end{equation}
where $F^2$ fits in an exact sequence
\begin{equation}
0 \to F_{n_1} \to F^2 \to F_{n_2} \to 0.
\end{equation}
Continueing this procedure, we get a sequence of stable complexes
\begin{equation}
{\Bbb E}={\Bbb E}^{(0)},{\Bbb E}^{(1)},\dots,{\Bbb E}^{(s)},\dots,
\end{equation} 
where ${\Bbb E}^{(s)}$ fits in an exact triangle
\begin{equation}
F^s \to {\Bbb E} \to {\Bbb E}^{(s)} \to F^s[1],
\end{equation}
$F^s \in {\cal S}$.
Since $v({\Bbb E}^{(s)})=v({\Bbb E}^{(0)})-\sum_i v(E_{n_i})$,
if $S$ generate a negative definite lattice
or $\langle \delta,v({\Bbb E}) \rangle>0$, then
$\langle v({\Bbb E}^{(s)})^2 \rangle <-(1+p_g)$
for some $s$. 
By Lemma \ref{lem:non-empty}, this is impossible.
Hence $\Hom(E_i,{\Bbb E}^{(s)})=0$, $1 \leq i \leq n$ 
for some $s$.

For a non-zero morphism $\psi:{\Bbb E} \to E_{n_0}[1]$,
we set ${\Bbb E}^{(-1)}[1]:=[{\Bbb E} \to E_{n_0}[1]]$.
Then ${\Bbb E}^{(-1)}$ fits in an exact triangle:
\begin{equation}
 E_{n_0} \to {\Bbb E}^{(-1)} \to {\Bbb E} \to E_{n_0}[1].
\end{equation}
\begin{claim}\label{claim:2}
${\Bbb E}^{(-1)}$ is a stable complex.
\end{claim}

Proof of Claim \ref{claim:2}:
For a non-zero morphism $\psi:{\Bbb E} \to E_{n_0}[1]$,
we have an exact sequence
\begin{equation}
\Hom({\Bbb E}[1],E_{n_0}[1]) \to
\Hom({\Bbb E}^{(-1)}[1],E_{n_0}[1]) \to {\Bbb C} 
\overset{\psi}{\to}
\Hom({\Bbb E},E_{n_0}[1]).
\end{equation}
By Lemma \ref{lem:hom[-1]},
$\Hom({\Bbb E}[1],E_{n_0}[1])=0$, and hence
$\Hom({\Bbb E}^{(-1)}[1],E_{n_0}[1])=0$.
By our assumption, ${\Bbb E}^{(-1)}$ satisfies \eqref{eq:length-cond}
and $H^i({\Bbb E}^{(-1)})$, $i=-1,0$ fits in the exact sequence
\begin{equation}
 \begin{CD}
  0 @>>> H^{-1}({\Bbb E}^{(-1)}) @>>> H^{-1}({\Bbb E}) @>>> @.\\
  E_{n_0} @>>> H^0({\Bbb E}^{(-1)}) @>>> H^0({\Bbb E}) @>>> 0.    
 \end{CD}
\end{equation}
If $H^{-1}({\Bbb E})=0$, then since
$\Hom(H^0({\Bbb E}^{(-1)}),E_{n_0})=\Hom({\Bbb E}^{(-1)}[1],E_{n_0}[1])=0$,
Lemma \ref{lem:key1} (2) implies that $H^0({\Bbb E}^{(-1)})$ is stable.
Assume that $H^{-1}({\Bbb E})\ne 0$.
If $H^{-1}({\Bbb E}) \to E_{n_0}$ is a zero map, then
since $H^0({\Bbb E})$ is of 0-dimensional and $E_{n_0}$ is 
locally free, we get
$\Ext^1(H^0({\Bbb E}),E_{n_0})=0$.
Hence the second line splits, which is a contradiction.
Thus $\xi:H^{-1}({\Bbb E}) \to E_{n_0}$ is non trivial.
Then by applying Lemma \ref{lem:key1} (3) to
$\xi^{\vee}:E_{n_0}^{\vee} \to H^{-1}({\Bbb E})^{\vee}$,
we see that
(1) $\xi^{\vee}$ is injective except finite subset of $X$
and $\coker(\xi^{\vee})$ is $\mu$-stable torsion free sheaf,
or (2)
$\xi^{\vee}$ is injective except a divisor of $X$
and $\coker(\xi^{\vee})$ is $\mu$-stable purely 1-dimensional sheaf,
or (3) $\xi^{\vee}$ is surjective in codimension 1
and $\ker \xi^{\vee}$ is a $\mu$-stable sheaf.
In the case of (1),
$H^0({\Bbb E}^{(-1)})$ is 0-dimensional and
$H^{-1}({\Bbb E}^{(-1)})$ is a $\mu$-stable sheaf.
If the case (2) occur, then 
$H^0({\Bbb E}^{(-1)})$ is a $\mu$-stable 1-dimensional sheaf
and $H^{-1}({\Bbb E}^{(-1)})=0$.
In the last case,
$H^0({\Bbb E}^{(-1)})$ is a $\mu$-stable torsion free sheaf and
$H^{-1}({\Bbb E}^{(-1)})=0$.
Therefore ${\Bbb E}^{(-1)}$ is a stable complex
and we complete the proof of the claim.

If there is a non-zero homomorphism
${\Bbb E}^{(-1)} \to E_{n_{-1}}[1]$, we set 
${\Bbb E}^{(-2)}[1]:=[{\Bbb E}^{(-1)} \to E_{n_{-1}}[1]]$.
Continueing this procedure, we get a sequence of stable complexes
\begin{equation}
 \dots, {\Bbb E}^{(-t)},\dots,{\Bbb E}^{(-1)},{\Bbb E}^{(0)}.
\end{equation} 
Since $v({\Bbb E}^{(-t)})=v({\Bbb E}^{(0)})+\sum_i v(E_{n_i})$,
if $S$ generate a negative definite lattice
or $\langle \delta,v({\Bbb E}) \rangle<0$,
then we see that $\Hom(E_i,{\Bbb E}^{(-t)})=0$, $1 \leq i \leq n$ 
for some $-t$. 
Therefore Lemma \ref{lem:univ-exist} holds.
\end{proof}

\begin{lem}\label{lem:obstruction}
Assume that $S$ satisfies the condition of
(i) or (ii) in Lemma \ref{lem:univ-exist}.
If there is an exact triangle
\begin{equation}
F \to {\Bbb E} \to {\Bbb E}' \to F[1],
\end{equation}
where $F \in {\cal S}$. 
Then 
$\Hom({\Bbb E},{\Bbb E}')\cong{\Bbb C}$ and
$\Hom({\Bbb E},{\Bbb E}' \otimes K_X) \cong H^0(X,K_X)$.
\end{lem}

\begin{proof}
We only show the first assertion.
We assume that there is a universal division
$\widetilde{\Bbb E}$.
Since $\Hom({\Bbb E}',\widetilde{\Bbb E}) \to 
\Hom({\Bbb E},\widetilde{\Bbb E})$
is surjective,
we have an exact triangle
\begin{equation}
F' \to {\Bbb E}' \to \widetilde{\Bbb E} \to  F'[1],
\end{equation}
where $F' \in {\cal S}$.
By the exact sequence
\begin{equation}
\Hom({\Bbb E},F') \to 
\Hom({\Bbb E},{\Bbb E}')
\to \Hom({\Bbb E},\widetilde{\Bbb E})={\Bbb C}
\end{equation}
and Lemma \ref{lem:hom[-1]},
we get our claim.
If there is a universal extension $\widehat{\Bbb E}$, 
we also see that $\Hom({\Bbb E},{\Bbb E}') \cong {\Bbb C}$.
\end{proof}

\subsection{Coherent systems}

We set
\begin{equation}
{\frak P}_{E_i}^{(n)}(v):=\{({\Bbb E},U)|{\Bbb E} \in M_H(v), U  \subset 
\Hom(E_i,{\Bbb E}), \dim U=n \}.
\end{equation}

${\frak P}_{E_i}^{(n)}(v)$ is the moduli space of coherent systems.
For the construction of ${\frak P}_{E_i}^{(n)}(v)$, 
see section \ref{subsect:coh-system}.
The Zariski tangent space of ${\frak P}_{E_i}^{(n)}(v)$ at $({\Bbb E},U)$ is
\begin{equation}
\Ext^1(U \otimes E_i \to {\Bbb E},{\Bbb E})/\End(U \otimes E_i)
\end{equation} 
and the obstruction for the infinitesimal deformation
belongs to the kernel of 
\begin{equation}
\tau:\Ext^2(U \otimes E_i \to {\Bbb E},{\Bbb E}) \to \Ext^2({\Bbb E},{\Bbb E})
\overset{\tr}{\to} H^2(X,{\cal O}_X).
\end{equation}
By Lemma \ref{lem:obstruction} and the Serre duality,
$\ker \tau=0$. Thus ${\frak P}_{E_i}^{(n)}(v)$ is a smooth
scheme with 
\begin{equation}
\begin{split}
\dim {\frak P}_{E_i}^{(n)}(v)&=
\dim \Ext^1(U \otimes E_i \to {\Bbb E},{\Bbb E})/\End(U \otimes E_i)\\
&=\langle v-nv_i,v \rangle-n^2+(1+p_g)\\
&=\frac{1}{2}(\dim M_H(v)+\dim M_H(v-nv_i)).
\end{split}
\end{equation}
For $({\Bbb E},U) \in {\frak P}_{E_i}^{(n)}(v)$, 
${\Bbb E}$ and
$[U \otimes E_i \to {\Bbb E}]$ are stable.
Hence we have morphisms
$\pi:{\frak P}_{E_i}^{(n)}(v) \to M_H(v)$ and
$\varpi:{\frak P}_{E_i}^{(n)}(v) \to M_H(v-nv_i)$.
\begin{rem}
We set ${\Bbb F}:=[U \otimes E_i \to {\Bbb E}]$.
Since $\Hom({\Bbb E}[1],E_i[1])=\Hom(E_i \otimes U,E_i[1])=0$, 
by the exact triangle
\begin{equation}
U \otimes E_i \to {\Bbb E} \to {\Bbb F} \to U \otimes E_i[1],
\end{equation}
we have an exact sequence
\begin{equation}
0 \to U^{\vee} \to \Hom({\Bbb F},E_i[1]) \to
\Hom({\Bbb E},E_i[1]) \to 0.
\end{equation}
Thus we have
\begin{equation}\label{eq:dual-system}
{\frak P}_{E_i}^{(n)}(v):=
\{({\Bbb F},U^{\vee})| {\Bbb F} \in M_H(v-nv_i), 
U^{\vee}  \subset 
\Hom({\Bbb F},E_i[1]), \dim U=n \}.
\end{equation}
\end{rem}
We set $F_i:=U \otimes E_i$.
Then we have the following exact and commutative diagram:
\begin{equation}\label{eq:tangent-map}
\begin{CD}
@. \Hom(F_i,F_i) @. @. @.\\
@. @VVV @.@.\\
\Hom(F_i \to {\Bbb E},F_i \to {\Bbb E}) @>>> 
\Ext^1(F_i \to {\Bbb E},F_i) @>>> \Ext^1(F_i \to {\Bbb E},{\Bbb E}) @>>>
\Ext^1(F_i \to {\Bbb E},F_i \to {\Bbb E})\\
@VVV @VVV @VVV @VVV \\
\Hom({\Bbb E},F_i \to {\Bbb E}) @>>> \Ext^1({\Bbb E},F_i) @>>> \Ext^1({\Bbb E},{\Bbb E}) @>>>
\Ext^1({\Bbb E},F_i \to {\Bbb E})
\end{CD}
\end{equation}
By Lemma \ref{lem:obstruction}, we see that 
$\Ext^1({\Bbb E},F_i) \to \Ext^1({\Bbb E},{\Bbb E})$ is injective,
which implies that 
\begin{equation}
\Ext^1(F_i \to {\Bbb E},{\Bbb E})/\End(F_i) \to \Ext^1(F_i \to {\Bbb E},F_i \to {\Bbb E})
\oplus \Ext^1({\Bbb E},{\Bbb E})
\end{equation}
is injective.
Therefore $\pi \times \varpi:
{\frak P}_{E_i}^{(n)}(v) \to M_H(v) \times M_H(v-nv_i)$ 
is a closed immersion.

\begin{defn}
\begin{equation}
M_H(v)_{E_i,n}:=\{{\Bbb E} \in M_H(v)| 
\dim \Hom(E_i,{\Bbb E})=n \}.
\end{equation}
\end{defn}
Then $\pi_*({\frak P}_{E_i}^{(n)}(v))=\cup_{k \geq n} M_H(v)_{E_i,n}$.

\section{An action of a Lie algebra}\label{sect:action}

We define a lattice  
\begin{equation}
L(S):=\left(\sum_{i=1}^n{\Bbb Z}v(E_i),-\langle \quad,\quad \rangle \right).
\end{equation}
Let ${\frak g}$ the Lie algebra associated to $L(S)$,
that is, the Cartan matrix of ${\frak g}$ is
$(-\langle v(E_i),v(E_j) \rangle_{i,j=1}^n)$.
In the same way as in \cite{Na:1998} and \cite{N:2002},
we shall construct an action of ${\frak g}$
on $\bigoplus_v H_*(M_H(v),{\Bbb C})$, where $v$ runs a suitable set of
Mukai vectors with \eqref{eq:min-cond}.

The fundamental class of ${\frak P}_{E_i}^{(n)}$ defines an operator
$f_{v_i}^{(n)}$:
\begin{equation}
\begin{matrix}
H_*(M_H(v-nv_i),{\Bbb C})& \to & H_*(M_H(v),{\Bbb C})\\
x & \mapsto & p_{2*}(p_1^*(x) \cap [{\frak P}_{E_i}^{(n)}(v)])
\end{matrix}
\end{equation}
where $p_1,p_2$ are the first and the second projections 
of $M_H(v-nv_i) \times M_H(v)$.
We also define the operator $e_{v_i}^{(n)}$:
\begin{equation}
\begin{matrix}
H_*(M_H(v),{\Bbb C})& \to & H_*(M_H(v-nv_i),{\Bbb C})\\
x & \mapsto & (-1)^{nr(v)}p_{1*}(p_2^*(x) \cap [{\frak P}_{E_i}^{(n)}(v)])
\end{matrix}
\end{equation}
where $r(v)=\frac{1}{2}(\dim M_H(v-v_i)-\dim M_H(v))=
-\langle v_i,v \rangle-1$.
We set $e_{v_i}:=e_{v_i}^{(1)}$ and $f_{v_i}:=f_{v_i}^{(1)}$.
We also set 
\begin{equation}
h_{v_i|H_*(M_H(v),{\Bbb C})}=\langle v_i,v \rangle \id_{H_*(M_H(v),{\Bbb C})}.
\end{equation}

\begin{thm}\label{thm:action}
Assume that $S$ satsifies the assumptions in (i) or (ii) of 
Lemma \ref{lem:univ-exist}.
Then $e_{v_i},f_{v_j},h_{v_k}$ satisfy the following relations:
%
\begin{align}
\label{eq:relation1}
[h_{v_i},e_{v_j}]&=-\langle v_i,v_j \rangle e_{v_j}\\
\label{eq:relation2}
[h_{v_i},f_{v_j}]&=\langle v_i,v_j \rangle f_{v_j}\\
\label{eq:relation3}
[e_{v_i},f_{v_j}]&=\delta_{i,j} h_{v_i},\\
\label{eq:serre}
\ad(e_{v_i})^{1+\langle v_i, v_j \rangle}(e_{v_j})& =
\ad(f_{v_i})^{1+\langle v_i, v_j \rangle}(f_{v_j})=0,\;i \ne j,
\end{align}
where $\ad$ means the adjoint action
$\ad(x)(y):=[x,y]=xy-yx$.
\end{thm}
Since $\langle (v\pm n v_i)^2 \rangle<-(1+p_g)$ for $n \gg0$,
$e_{v_i}$ and $f_{v_i}$ are locally nilpotent.
Therefore we get an integral representation of ${\frak g}$.

\subsection{Proof of Theorem \ref{thm:action}}
The proof is similar to \cite{Na:1998} and
\cite{N:2002}.
We first note that the Serre relations \eqref{eq:serre} 
follows from the other relations and 
$\langle (v\pm n v_i)^2 \rangle<-(1+p_g)$ for $n \gg0$:
Let $L$ be the subspace 
of $\Hom(\bigoplus_{k \in {\Bbb Z}}H_*(M_H(v+kv_i)),
\bigoplus_{k \in {\Bbb Z}}H_*(M_H(v+v_j+kv_i)))$ generated by
$\ad(e_{v_i})^n (e_{v_j})$, $n \geq 0$.
Then $\frak{sl}_2$ generated by $e_{v_i},f_{v_i},h_{v_i}$
acts on $L$.
Since $\langle (v\pm n v_i)^2 \rangle<-(1+p_g)$ for $n \gg0$,
$L$ is of finite dimension.
By the theory of the $\frak{sl}_2$-representation,
we get 
$\ad(e_{v_i})^{1+\langle v_i, v_j \rangle}(e_{v_j})$.
The proof of the other relation is the same.

Hence we only need to show relations \eqref{eq:relation1},
\eqref{eq:relation2} and \eqref{eq:relation3}.
The proof of \eqref{eq:relation1}, \eqref{eq:relation2} are easy.
We shall prove \eqref{eq:relation3}.
For this purpose, we shall study the convolution products:
\begin{equation}\label{eq:products}
\begin{gathered}
   p_{13*}\left(p_{12}^* \left[\omega({\frak P}_{E_i}^{(n_i)}(v))\right]\cap
     p_{23}^* \left[{\frak P}_{E_j}^{(n_j)}(v)\right]
   \right),
\\
   q_{13*}\left(q_{12}^* \left[{\frak P}_{E_j}^{(n_j)}(v-n_i v_i)\right]\cap
     q_{23}^* \left[\omega({\frak P}_{E_i}^{(n_i)}(v-n_i v_i))\right]
   \right),
\end{gathered}
\end{equation}
where $p_{ij}$ and $q_{ij}$ are projections to the product of $i$-th
and $j$-th factors in
\begin{equation*}
   M_H(v-n_i v_i)\times M_H(v) \times M_H(v-n_j v_j),
\qquad
   M_H(v-n_i v_i)\times M_H(v-n_i v_i-n_j v_j) \times M_H(v-n_j v_j)
\end{equation*}
respectively, and $\omega$ is the exchange of the factor. The both
products have degree
$\frac{1}{2}(\dim M_H(v-n_i v_i)+\dim M_H(v-n_j v_j)).
$

(I) We first study the case where $i \ne j$.
\begin{lem}\label{lem:commutative1}
We have an isomophism over 
$M_H(v-n_i v_i) \times M_H(v-n_j v_j)$:
\begin{equation}
p_{12}^{-1}(\omega({\frak P}_{E_i}^{(n_i)}(v))) \cap 
p_{23}^{-1}({\frak P}_{E_j}^{(n_j)}(v))
\to 
q_{12}^{-1}({\frak P}_{E_j}^{(n_j)}(v-n_iv_i)) \cap 
q_{23}^{-1}(\omega({\frak P}_{E_i}^{(n_i)}(v-n_j v_j))).
\end{equation}
\end{lem}

\begin{proof}
For $F_i \to {\Bbb E}_1$ and $[F_i \to {\Bbb E}_1] \to F_j[1]$,
we set 
\begin{equation}
\begin{split}
{\Bbb E}_2:=& [[F_i \to {\Bbb E}_1] \to F_j[1]][-1],\\
{\Bbb E}:=& [{\Bbb E}_1 \to F_j[1]][-1].
\end{split}
\end{equation}
Applying the Octahedral axiom to
${\Bbb E}_1 \to [F_i \to {\Bbb E}_1] \to F_j[1]$,
 we have a commutative diagram of exact triangles:
\begin{equation}\label{eq:comm}
\begin{CD}
@. F_i @= F_i @. @. \\
@. @VVV @VVV @. @. \\
F_j @>>> {\Bbb E} @>>> {\Bbb E}_1 @>>>
F_j[1] \\
@| @VVV @VVV @|\\
F_j @>>> {\Bbb E}_2 @>>> [F_i \to {\Bbb E}_1] @>>> F_j[1]\\
@. @VVV @VVV @. @. \\
@. F_i[1] @= F_i[1] @. @. 
\end{CD}
\end{equation}
Hence ${\Bbb E}_1 \cong [F_j \to {\Bbb E}]$,
${\Bbb E}_2 \cong [F_i \to {\Bbb E}]$ and 
$[F_i \to {\Bbb E}_1] \cong [F_i \oplus F_j \to {\Bbb E}]$.
Conversely for ${\Bbb E}:=[{\Bbb E}_1 \to F_j[1]][-1]$ and
${\Bbb E}_2:=[F_i \to {\Bbb E}]$, we get the commutative 
diagram of exact triangles
\eqref{eq:comm}.
Since the correspondence is functorial, we have 
a desired isomorphism
\begin{equation}
p_{12}^{-1}(\omega({\frak P}_{E_i}^{(n_i)}(v))) \cap 
p_{23}^{-1}({\frak P}_{E_j}^{(n_j)}(v))
\to 
q_{12}^{-1}({\frak P}_{E_j}^{(n_j)}(v-n_iv_i)) \cap 
q_{23}^{-1}(\omega({\frak P}_{E_i}^{(n_i)}(v-n_j v_j))).
\end{equation}
\end{proof}

\begin{lem}\label{lem:commutative2}
\begin{enumerate}
\item
$p_{12}^{-1}(\omega({\frak P}_{E_i}^{(n_i)}(v))) \cap 
p_{23}^{-1}({\frak P}_{E_j}^{(n_j)}(v))
\to M_H(v-n_i v_i) \times M_H(v-n_j v_j)
$
is injective.
\item
$q_{12}^{-1}({\frak P}_{E_j}^{(n_j)}(v-n_iv_i)) \cap 
q_{23}^{-1}(\omega({\frak P}_{E_i}^{(n_i)}(v-n_j v_j)))
\to M_H(v-n_i v_i) \times M_H(v-n_j v_j)$
is injective.
\end{enumerate}
\end{lem}

\begin{proof}
We shall prove (i).
The proof of (ii) is similar.
Assume that we have isomorphisms in the derived category:
\begin{equation}
\begin{split}
[F_1 \to {\Bbb E}] & \cong [F_1 \to {\Bbb E}'],\\
[F_2 \to {\Bbb E}] & \cong [F_2 \to {\Bbb E}'].
\end{split}
\end{equation}
We shall show that there is an isomorphism 
$\phi:{\Bbb E} \to {\Bbb E}'$ which is compatible with the morphisms
$F_i \to {\Bbb E}$, $F_i \to {\Bbb E}'$. 
Applying $\Hom({\Bbb E}',\quad)$ to the exact triangles
\begin{equation}
\begin{CD}
F_1 @>>> {\Bbb E} @>>> [F_1 \to {\Bbb E}] @>>> F_1[1]\\
@| @VVV @VVV @|\\
F_1 @>>> [F_2 \to {\Bbb E}] @>>> [F_1 \oplus F_2 \to {\Bbb E}] @>>> F_1[1],
\end{CD}
\end{equation}
we get a commutative diagram
\begin{equation}
\begin{CD}
\Hom({\Bbb E}',F_1) @>>> \Hom({\Bbb E}',{\Bbb E}) @>>> 
\Hom({\Bbb E}',F_1 \to {\Bbb E}) @>>> \Ext^1({\Bbb E}',F_1)\\
@| @VVV @VVV @|\\
\Hom({\Bbb E}',F_1) @>>> \Hom({\Bbb E}',F_2 \to {\Bbb E}) @>>> 
\Hom({\Bbb E}', F_1 \oplus F_2 \to {\Bbb E}) @>>> 
\Ext^1({\Bbb E}',F_1).
\end{CD}
\end{equation}
Since $\Hom({\Bbb E}',F_1)=0$, Lemma \ref{lem:obstruction}
implies that
\begin{equation}
\Hom({\Bbb E}',F_2 \to {\Bbb E}) \to 
\Hom({\Bbb E}',F_1 \oplus F_2 \to {\Bbb E})
\end{equation} 
is an isomorphism.
Hence $\Hom({\Bbb E}',{\Bbb E}) \to \Hom({\Bbb E}',F_1 \to {\Bbb E}) 
\cong {\Bbb C}$ 
is also an isomorphism.
We also have an isomorphism 
\begin{equation}
\Hom({\Bbb E},{\Bbb E}') \to 
\Hom({\Bbb E},F_1 \to {\Bbb E}').
\end{equation}
Then the claim easily follow from these isomorphisms.
Hence 
\begin{equation}
p_{12}^{-1}(\omega({\frak P}_{E_i}^{(n_i)}(v))) \cap 
p_{23}^{-1}({\frak P}_{E_j}^{(n_j)}(v))
\to M_H(v-n_i v_i) \times M_H(v-n_j v_j)
\end{equation}
is injective.
\end{proof}

\begin{lem}\label{lem:commutative3}
If $i \ne j$, then $p_{12}^{-1}(\omega({\frak P}_{E_i}^{(n_i)}(v)))$ and 
$p_{23}^{-1}({\frak P}_{E_j}^{(n_j)}(v))$ 
intersect transversely and
$q_{12}^{-1}({\frak P}_{E_j}^{(n_j)}(v-n_iv_i))$ and 
$q_{23}^{-1}(\omega({\frak P}_{E_i}^{(n_i)}(v-n_j v_j)))$ intersect 
transversely.
\end{lem}

\begin{proof}
We set $F_i:=U_i \otimes E_i$ and $F_j:=U_j \otimes E_j$.
We shall show that the map of the tangent spaces 
\begin{equation}\label{eq:trans1}
\begin{split}
& \Ext^1(F_i \to {\Bbb E},{\Bbb E})/\End(F_i) \oplus
\Ext^1(F_j \to {\Bbb E},{\Bbb E})/\End(F_j)\\
 \to & 
\Ext^1({\Bbb E},{\Bbb E})
\end{split}
\end{equation}
is surjective and 
\begin{equation}\label{eq:trans2}
\begin{split}
& \Ext^1(F_i \oplus F_j \to {\Bbb E},F_j \to {\Bbb E})/\End(F_i) \oplus
\Ext^1(F_i \oplus F_j \to {\Bbb E},F_i \to {\Bbb E})/\End(F_j)\\
 \to & 
\Ext^1(F_i \oplus F_j \to {\Bbb E},F_i \oplus F_j \to {\Bbb E})
\end{split}
\end{equation}
is surjective.

We shall only prove \eqref{eq:trans1}.
By \eqref{eq:tangent-map} and Lemma \ref{lem:obstruction},
it is sufficient to show that the natural homomorphism
\begin{equation}
\Ext^1(F_i \to {\Bbb E},{\Bbb E}) \to
\Ext^1({\Bbb E},{\Bbb E}) \to \Ext^1(F_j,{\Bbb E})
\end{equation}
is surjective.
Since $\Ext^2(F_i \oplus F_j \to {\Bbb E},{\Bbb E}) 
\cong \Ext^2(F_i \to {\Bbb E},{\Bbb E}) \cong H^2(X,{\cal O}_X)$,
the exact triangle
\begin{equation}
F_j \to [F_i \to {\Bbb E}] \to [F_i \oplus F_j \to {\Bbb E}] \to 
F_j[1]
\end{equation}
implies that 
this homomorphism is surjective. 
\end{proof}
By Lemmas \ref{lem:commutative1}, \ref{lem:commutative2},
\ref{lem:commutative3}, we obtain that

\begin{equation}\label{eq:products2}
   p_{13*}\left(p_{12}^* \left[\omega({\frak P}_{E_i}^{(n_i)}(v))\right]\cap
     p_{23}^* \left[{\frak P}_{E_j}^{(n_j)}(v)\right]
   \right)=
   q_{13*}\left(q_{12}^* \left[{\frak P}_{E_j}^{(n_j)}(v-n_i v_i)\right]\cap
     q_{23}^* \left[\omega{\frak P}_{E_i}^{(n_i)}(v-n_i v_i)\right]
   \right).
\end{equation}
Hence we get 
\begin{equation}
[e_{v_i},f_{v_j}]=0,\; i \ne j.
\end{equation}

(II) We next treat the case where $i=j$.
This case was treated by Nakajima \cite{N:2002}.
For convenience of the reader, we write a self-contained proof.
We assume that $n=1$.
If $ i=j$, then 
$p_{12}^{-1}(\omega({\frak P}_{E_i}^{(1)}))$ and 
$p_{23}^{-1}({\frak P}_{E_j}^{(1)})$ 
intersect transversely outside $p_{13}^{-1}(\Delta_{M_H(v-v_i)})$,
and 
$q_{12}^{-1}({\frak P}_{E_i}^{(1)}(v-v_i))$ and
$q_{23}^{-1}(\omega{\frak P}_{E_i}^{(1)}(v-v_i))$
intersect transversely outside $q_{13}^{-1}(\Delta_{M_H(v-v_i)})$.
Then we see that
\begin{equation}\label{eq:products3}
   p_{13*}\left(p_{12}^* \left[\omega({\frak P}_{E_i}^{(1)}(v))\right]\cap
     p_{23}^* \left[{\frak P}_{E_i}^{(1)}(v)\right]
   \right)=
   q_{13*}\left(q_{12}^* \left[{\frak P}_{E_i}^{(1)}(v-v_i)\right]\cap
     q_{23}^* \left[\omega{\frak P}_{E_i}^{(1)}(v-v_i)\right]
   \right)
\end{equation}
outside $\Delta_{M_H(v-v_i)}$. 
Thus
\begin{equation}\label{eq:products4}
   p_{13*}\left(p_{12}^* \left[\omega({\frak P}_{E_i}^{(1)}(v))\right]\cap
     p_{23}^* \left[{\frak P}_{E_i}^{(1)}(v)\right]
   \right)=
   q_{13*}\left(q_{12}^* \left[{\frak P}_{E_i}^{(1)}(v-v_i)\right]\cap
     q_{23}^* \left[\omega{\frak P}_{E_i}^{(1)}(v-v_i)\right]
   \right)+c \Delta_{M_H(v-v_i)}
\end{equation}
for some integer $c$.
In order to compute $c$, we may restrict to a suitable open neighbourhood 
of the generic point of
$\Delta_{M_H(v-v_i)}$. We set $w:=v-v_i$.

(II-1)
Assume that $-\langle v_i,w \rangle \geq 0$.
We set
\begin{equation}
\begin{split}
M_H(w)':=& M_H(w)_{E_i,-\langle v_i,w \rangle},\\
M_H(w-v_i)':=& M_H(w-v_i) \setminus
\pi({\frak P}_{E_i}^{(-\langle v_i,w \rangle)}(w-v_i)).
\end{split}
\end{equation}
Then ${\frak P}_{E_i}^{(1)}(w)':=\pi^{-1}(M_H(w)')$
is a projective bundle over $M_H(w)'$
and ${\frak P}_{E_i}^{(1)}(w)' \to M_H(w-v_i)'$ is a closed immersion.
We have a fiber product diagram:
\begin{equation}
\begin{CD}
{\frak P}_{E_i}^{(1)}(w)' @>>> q^{-1}_{23}(\omega({\frak P}_{E_i}^{(1)}(w)'))\\
@VVV @VVV \\
q^{-1}_{12}({\frak P}_{E_i}^{(1)}(w)') @>>>
M_H(w)' \times M_H(w-v_i)' \times M_H(w)'.
\end{CD}
\end{equation}
By the excess intersection theory, we get that
\begin{equation}
q^*_{12}\left[{\frak P}_{E_i}^{(1)}(w)' \right]
\cap q^*_{23}\left[\omega({\frak P}_{E_i}^{(1)}(w)')\right]
=c_{top}(N_{{\frak P}_{E_i}^{(1)}(w)'/M_H(w-v_i)'}).
\end{equation} 
We take ${\Bbb E} \in M_H(w)$ with $\Ext^1(E_i,{\Bbb E})=0$.
We set $V:=\Hom(E_i,{\Bbb E})$.
Let ${\Bbb P}:={\Bbb P}(V^{\vee})$ be the fiber of $\pi$.
Then 
\begin{equation}
q_{13*}\left(q^*_{12}\left[{\frak P}_{E_i}^{(1)}(w)' \right]
\cap q^*_{23}\left[\omega({\frak P}_{E_i}^{(1)}(w)')\right] \right)
=\left(\int_{{\Bbb P}} 
c_{top}(N_{{\frak P}_{E_i}^{(1)}(w)'/M_H(w-v_i)'})\right)
\Delta_{M_H(w)'}.
\end{equation}
We have a family of non-trivial homomorphisms:
\begin{equation}
{\cal O}_{\Bbb P}(-1) \boxtimes E_i \to {\cal O}_{\Bbb P} \boxtimes {\Bbb E} .
\end{equation}
We set 
\begin{equation}
 {\cal E}:=
 [{\cal O}_{\Bbb P}(-1) \boxtimes E_i \to 
 {\cal O}_{\Bbb P} \boxtimes {\Bbb E}].
\end{equation}
We have an exact sequence
\begin{equation}
 \begin{CD}
  \Ext^1_{p_{\Bbb P}}({\cal E},{\cal O}_{\Bbb P} \boxtimes{\Bbb E}) @>>> 
  \Ext^1_{p_{\Bbb P}}({\cal E},{\cal E}) @>>>
  \Ext^2_{p_{\Bbb P}}({\cal E},{\cal O}_{\Bbb P}(-1) \boxtimes E_i)\\
  \Ext^2_{p_{\Bbb P}}({\cal E},{\cal O}_{\Bbb P} \boxtimes{\Bbb E}) @>>> 
  \Ext^2_{p_{\Bbb P}}({\cal E},{\cal E}). @.
 \end{CD}
\end{equation}
The restriction of the 
normal bundle $(N_{{\frak P}_{E_i}^{(1)}(w)'/M_H(w-v_i)'})_{|{\Bbb P}}$ is
\begin{equation}
 \Ext^2_{p_{\Bbb P}}({\cal E},{\cal O}_{\Bbb P}(-1) \boxtimes E_i)
 =\Hom_{p_{\Bbb P}}({\cal O}_{\Bbb P}(-1) \boxtimes E_i,{\cal E})^{\vee}.
 \end{equation}
By the exact triangle
\begin{equation}
 {\cal O}_{\Bbb P}(-1) \boxtimes E_i \to {\cal O}_{\Bbb P} \boxtimes {\Bbb E}
 \to {\cal E} \to {\cal O}_{\Bbb P}(-1) \boxtimes E_i[1],
 \end{equation}
we get an exact sequence
 \begin{equation}
0 \to {\cal O}_{\Bbb P} \to V \otimes {\cal O}_{\Bbb P}(1) \to 
\Hom_{p_{\Bbb P}}({\cal O}_{\Bbb P}(-1) \boxtimes E_i,{\cal E}) \to 0.
 \end{equation}
Hence 
$\Hom_{p_{\Bbb P}}({\cal O}_{\Bbb P}(-1) \boxtimes E_i,{\cal E})^{\vee}
=\Omega_{\Bbb P}^1$.
Therefore
\begin{equation}
\int_{\Bbb P} c_{top}(N_{{\frak P}_{E_i}^{(1)}(w)'/M_H(w-v_i)'})
=(-1)^{\dim {\Bbb P}}(\dim {\Bbb P}+1)
=(-1)^{-\langle v_i,w \rangle}\langle v_i,w \rangle.
\end{equation}  
Since ${\frak P}_{E_i}^{(1)}(v)'$ does not meet
$M_H(v) \times M_H(w)'$,
\begin{equation}
   p_{13*}\left(p_{12}^* \left[\omega({\frak P}_{E_i}^{(1)}(v))\right]\cap
     p_{23}^* \left[{\frak P}_{E_i}^{(1)}(v)\right]
   \right)=0
   \end{equation} 
   on $M_H(w)' \times M_H(w)'$.
  Hence we see that 
  \begin{equation}\label{eq:sl2}
  [e_{v_i},f_{v_i}]_{|H_*(M_H(w),{\Bbb C})}=
\langle v_i,w \rangle \id_{H_*(M_H(w),{\Bbb C})}=
  h_{v_i|H_*(M_H(w),{\Bbb C})}.
  \end{equation}

(II-2)  
Assume that $\langle v_i,w \rangle \geq 0$.
We set $M_H(v)':=M_H(v) \setminus \pi({\frak P}_{E_i}^{(2)})$ and
$M_H(w)':=M_H(w)_{E_i,0}$.
For ${\Bbb E} \in M_H(v)'$,
we set $V:=\Ext^1({\Bbb E},E_i)$.
We have a family of exact triangles:
\begin{equation}
{\cal O}_{\Bbb P} \boxtimes E_i \to {\cal E}' \to 
{\cal O}_{\Bbb P}(-1) \boxtimes {\Bbb E}
\to {\cal O}_{\Bbb P} \boxtimes E_i[1].
\end{equation}
The restriction of the normal bundle 
$(N_{{\frak P}_{E_i}^{(1)}(w)'/M_H(v)'})_{|{\Bbb P}}$ is 
\begin{equation}
\Ext^1_{p_{\Bbb P}}({\cal O}_{\Bbb P} \boxtimes E_i,{\cal E}')=
\Ext^1_{p_{\Bbb P}}({\cal E}',{\cal O}_{\Bbb P} \boxtimes E_i)^{\vee}.
\end{equation}
We have an exact sequence
\begin{equation}
0 \to {\cal O}_{\Bbb P}=
\Hom_{p_{\Bbb P}}({\cal O}_{\Bbb P} \boxtimes E_i,
{\cal O}_{\Bbb P} \boxtimes E_i) \to
\Ext^1_{p_{\Bbb P}}({\cal O}_{\Bbb P}(-1) \boxtimes {\Bbb E},
{\cal O}_{\Bbb P} \boxtimes E_i)
\to  \Ext^1_{p_{\Bbb P}}({\cal E}',{\cal O}_{\Bbb P} \boxtimes E_i)
\to 0.
\end{equation}
Hence $\Ext^1_{p_{\Bbb P}}({\cal O}_{\Bbb P} \boxtimes E_i,{\cal E}')=
\Omega_{\Bbb P}^1$.
 Therefore
\begin{equation}
\int_{\Bbb P} c_{top}(N_{{\frak P}_{E_i}^{(1)}(w)'/M_H(v)'})
=(-1)^{\dim {\Bbb P}}(\dim {\Bbb P}+1)
=-(-1)^{-\langle v_i,w \rangle}\langle v_i,w \rangle.
\end{equation}  
By using this equality, we see that \eqref{eq:sl2} also holds.

\subsection{The case where the twisted degree is zero}
\label{subsect:G-twisted}

Let $G$ be an element of $K(X)$.
\begin{defn}
Let ${\Bbb E} \in {\bf D}(X)$ be an object such that 
$\deg_G({\Bbb E})=0$ and
\begin{equation}
\chi_G({\Bbb E})=\min\{\chi_G(E')>0 |E' \in \Coh(X),\deg_G(E')=0 \}.
\end{equation}
${\Bbb E}$ is $G$-twisted stable, if
\begin{enumerate}
\item
$H^i({\Bbb E})=0$, $i \ne 0$ and $H^0({\Bbb E})$ is $G$-twisted stable, or 
\item
$H^i({\Bbb E})=0$, $i \ne -1$ and $H^{-1}({\Bbb E})$ 
is $G$-twisted stable.
\end{enumerate}

Let $M_H^G(v)$ be the moduli space of $G$-twisted stable complex
${\Bbb E}$ with $v({\Bbb E})=v$.
\end{defn}
\begin{rem}
If (ii) holds, then
\begin{equation}
\chi(H^{-1}({\Bbb E}))=\max
\{\chi_G(E')<0 |E' \in \Coh(X),\deg_G(E')=0 \}.
\end{equation}
\end{rem}
Let $E_i$, $i=1,\dots,n$ be a collection of $G$-twisted stable
vector bundles with $\deg_G(E_i)=\chi_G(E_i)=0$ and 
$\langle v(E_i)^2 \rangle=-2$.
Assume that $E_i$ satisfies the condition 
\eqref{eq:(-2)-condition1}. 
By using Lemmas \ref{lem:key2} and \ref{lem:key2'},
 we also obtain the same assertions in Lemma \ref{lem:univ-exist}.
Hence we also get an action of the Lie algebra associated to
$E_i, i=1,\dots,n$.

\section{Examples}\label{sect:example}

\subsection{Stable sheaves on a $K3$ surface}

Let $X$ be a $K3$ surface and $H$ an ample divisor on $X$.
Let $G$ be a semi-stable vector bundle with respect to $H$
such that $\langle v(G)^2 \rangle=0$.
Assume that $G=\bigoplus_{i=0}^n E_i^{\oplus a_i}$,
where $E_i$ is a $G$-twisted stable vector bundle such that
\begin{equation}
\begin{split}
\frac{\deg(E_i)}{\rk E_i}&=\frac{\deg(G)}{\rk G},\\
\frac{\chi_G(E_i)}{\rk E_i}&=\frac{\chi_G(G)}{\rk G}=0.
\end{split}
\end{equation}
By \cite[Thm. 0.1]{O-Y:1}, $v(E_0),v(E_1),\dots,v(E_n)$
generate a lattice of affine type.
We may assume that $a_0=1$. 
We set
\begin{equation}
l:=\min\{\deg_G(E)>0| E \in \Coh(X) \}.
\end{equation}
We set $v_i:=v(E_i)$, $i=0,1,\dots,n$.
Let ${\frak g}$ be the affine Lie algebra associated with
$v_i$, $i=0,1,\dots,n$ 
and $\overline{\frak g}$ the finite Lie algebra 
associated with $v_i$, $i=1,\dots,n$.
Let $\overline{\frak h}$ be the Cartan subalgebra of $\overline{\frak g}$. 
For a root $\alpha$, $\overline{\frak g}_{\alpha}$ 
denotes the root space of $\alpha$.
$\theta:=\sum_{i=1}^n a_i v_i$ denotes
the highest root of $\overline{\frak g}$.
Then ${\frak g}$ has the following standard expression:
\begin{equation}
{\frak g}={\Bbb C}[t,t^{-1}] \otimes \overline{\frak g} 
\oplus {\Bbb C}c \oplus {\Bbb C}d
\end{equation}
where   
\begin{alignat}{4}
e_{v_i}= & \;1 \otimes \overline{e}_{v_i}, \quad \;
f_{v_i}= & 1 \otimes \overline{e}_{-v_i},\quad
h_{v_i}= & \;1 \otimes \overline{h}_{v_i}
\quad & 1 \leq i \leq n,\label{eq:chevalley1}\\
e_{v_0}= & \; t \otimes \overline{e}_{-\theta},\quad
f_{v_0}=& \;t^{-1} \otimes \overline{e}_{\theta},\quad\;
h_{v_0}= & -\sum_{i=1}^n a_i h_{v_i}+c,& \quad
\end{alignat}
%
%
%
$\overline{e}_{\alpha} \in \overline{\frak g}_{\alpha}$,
$\overline{h}_{v_i} \in {\frak h}$ and
\eqref{eq:chevalley1} are the Chevalley generator of 
$\overline{\frak g}$.
Hence we get 
\begin{equation}
c=\sum_{i=0}^n a_i h_{v_i}.
\end{equation}
The action of $d$ on $H_*(M_H(v),{\Bbb C})$
is defined as follows:
We take $w \in H^*(X,{\Bbb Q})$ such that
$\langle w,v(E_i) \rangle=\delta_{i,0}$, $i=0,1,\dots,n$
and set 
\begin{equation}
d_{|H_*(M_H(v),{\Bbb C})}:=\langle w,v \rangle 
\id_{H_*(M_H(v),{\Bbb C})}.
\end{equation}
Then we have a desired properties:
\begin{equation}
\begin{split}
[d,e_{v_i}]&=\delta_{i,0} e_{v_i},\\
[d,f_{v_i}]&=-\delta_{i,0} f_{v_i}.
\end{split}
\end{equation}

\begin{prop}
Let ${\frak g}$ be the affine Lie algebra 
associated to $E_0,E_1,\dots,E_n$.
Assume that $E_i$ are $\mu$-stable for all $i$.
Then we have an action of 
${\frak g}$ on $\bigoplus_v H_*(M_H(v),{\Bbb C})$
such that the center $c$ acts as a scalar multiplication
$\langle v,v(G) \rangle$,
where $v$ is a Mukai vector with $\deg_G(v)=l$. 
\end{prop}

\begin{ex}\label{ex:1}
Let $C:=(-a_{i,j})_{i,j=0}^n$ be a Cartan matrix of affine type and
$\delta:=(a_0,a_1,\dots,a_n)$, $a_i \in {\Bbb Z}_{>0}$ 
the primitive vector with $\delta C=0$. 
Let $(X,H)$ be a polarized $K3$ surface such that
\begin{enumerate}
\item
$\Pic(X)=\bigoplus_{i=0}^n {\Bbb Z} \xi_i$,
$(\xi_i,\xi_j)=-a_{i,j}+2ra$ and
\item
$H=\sum_{i=0}^n a_i \xi_i$.
\end{enumerate}
For an existence of $(X,H)$, see \cite[sect. 3]{O-Y:1}.
We set $v_i:=r+\xi_i+a \rho$.
Then 
\begin{enumerate}
\item
$(\langle v_i,v_j \rangle)_{i,j=0}^n=-C$,
\item
$\deg(v_i)=(\xi_i,H)=2ra(\sum_{i=0}^n a_i)$
and 
\item
$v:=\sum_i a_i v_i$ is a primitive isotropic Mukai vector.
\end{enumerate}
\begin{lem}
Let $E_i$ be a $v$-twisted stable vector bundle with respect to
$H$ with
$v(E_i)=v_i$.
Then $E_i$ is $\mu$-stable.
\end{lem}

\begin{proof}
For a coherent sheaf $F$, 
we set $c_1(F):=\sum_i x_i \xi_i$. 
Then $\deg(F)=(\sum_i x_i)2ra(\sum_i a_i)$.
Since $\rk E_i=r$ and $\deg E_i=2ra(\sum_i a_i)$ for all $i$,
if $\deg(F)/\rk F=\deg(E_i)/\rk E_i=2a(\sum_i a_i)$,
then $\rk F=(\sum_i x_i)r \geq \rk E_i$.
Therefore $E_i$ are $\mu$-stable.
\end{proof}
We set $w:=(r(\sum_i x_i)-1)+\sum_i x_i \xi_i+a \rho$,
$x_i \in {\Bbb Z}$.
Then 
\begin{equation}
\deg_G(w)=\min\{\deg_G(E)>0| E \in \Coh(X) \},
\end{equation}
where $G \in M_H(v)$.
Hence
we have an action of 
${\frak g}$ on $\bigoplus_w H_*(M_H(w),{\Bbb C})$,
where $w=\sum_i x_i v_i-1+a \rho$,
$x_i,a \in {\Bbb Z}$.
\end{ex}

\vspace{1pc}

Let $G$ be a vector bundle such that $\rk G=(H^2)$ and
$c_1(G)=H$.
For a Mukai vector $v:=(1+(D,H))-D+a \rho$, we get
\begin{equation}\label{eq:min}
\begin{split}
\deg_{G}(v)=&(H,H)(-D,H)-(1+(D,H))(-H,H)\\
=&(H,H)\\
=&\min \{\deg_{G}(E')>0|E' \in \Coh(X) \}.
\end{split}
\end{equation}
%
%
Let $C_1,C_2,\dots,C_n$ be irreducible $(-2)$-curves on $X$.
We set $v_i:=(C_i,H)-C_i$.
\begin{lem}\label{lem:mu-stable}
There is a stable vector bundle $E_i$ 
with $v(E_i)=v_i$.
Moreover if 
$H=n H'$ and $(C_i,H')<2(n-1)({H'}^2)$, then 
$E_i$ is $\mu$-stable.
\end{lem}

\begin{proof}
There is a semi-stable sheaf $E_i$ with $v(E_i)=v_i$.
We shall show that $E_i$ is stable.
Let $\bigoplus_{j=1}^s E_{i,j}$ be the Jordan-H\"{o}lder
grading of $E_i$ with respect to the Gieseker stability.
We set $v(E_{i,j}):=r_j-D_j+a_j \rho$.
Then $(D_j,H)/r_j=1$ and $a_j/r_j=0$, and hence $(D_j,H)>0$
and $(D_j^2)=-2$, which implies that $D_j$ is effective.
By our assumption on $C_i$, $s=1$. Thus $E_i$ is stable.
Assume that $H=nH'$ and $(C_i,H')<2(n-1)({H'}^2)$.
Let $\bigoplus_{j=1}^s E_{i,j}$ be the Jordan-H\"{o}lder
grading of $E_i$ with respect to the $\mu$-stability.
We set $v(E_{i,j}):=r_j-D_j+a_j \rho$.
Then $r_j=(D_j,H)=n(D_j,H')$, and hence $(D_j,H')>0$.
By the stability of $E_{i,j}$, $\langle v(E_{i,j})^2 \rangle
=(D_j^2)-2r_j a_j \geq -2$.
By the Hodge index theorem, $(D_j,H')^2 \geq (D_j^2)({H'}^2)$.
If $a_j >0$, then we see that
$(C_i,H')>(D_j,H') \geq 2(n-1)({H'}^2)$.
Therefore $a_j \leq 0$. Since $\sum_j a_j=0$,
$a_j=0$ for all $j$.
Since $E_i$ is stable, $s=1$. Thus $E_i$ is $\mu$-stable.
\end{proof}


\begin{prop} 
Assume that $E_i$ are $\mu$-stable.
Then we have an action of the Lie algebra ${\frak g}$
associated to $C_i$, $i=1,2,\dots,n$ on
$\bigoplus_v H_*(M_H(v),{\Bbb C})$, where $v=(1+(D,H))-D+a \rho$,
$D \in \Pic(X)$, $a \in {\Bbb Z}$.
\end{prop}

\begin{ex}\label{ex:2}
Let $\pi:X \to {\Bbb P}^1$ be an elliptic $K3$ surface with a section
$C_0$.
Let $C_1,\dots,C_n$ be smooth $(-2)$-curves on fibers of $\pi$.
We set $v_i:=(C_i,H)-C_i$, $i=0,1,\dots,n$.
Then $(\langle v_i,v_j \rangle_{i,j})=((C_i,C_j)_{i,j})$. 
We assume that $(C_i,C_j) \leq 1$.
Hence we get an action of the Lie algebra generated by 
$C_i, 0 \leq i \leq n$ on  
$\bigoplus_v H_*(M_H(v),{\Bbb C})$, 
where $v=(1+(D,H))-D+a \rho$, $D \in \Pic(X)$,
$a \in {\Bbb Z}$.
\end{ex}

\vspace{1pc}

\begin{ex}\label{ex:3}

In the notation of Example \ref{ex:1}, we set 
$v_i:=v(E_i)$.
Then we see that
\begin{equation}
\begin{split}
&\{\chi_G(w)|w \in v({\bf D}(X)), \deg_G(w)=0 \}\\
=&\{\chi_G(w)|w=\sum_i x_i v_i +y \rho \}\\
=&{\Bbb Z}\langle v(G),\rho \rangle.
\end{split}
\end{equation}
Hence we have an action of ${\frak g}$ on
$\bigoplus_w H_*(M_H^G(w),{\Bbb C})$, where
$w=\sum_i x_i v_i+\rho$, $x_i \in {\Bbb Z}$.

Let $\pi:X \to {\Bbb P}^1$ be the elliptic $K3$ surface 
as in Example \ref{ex:2}.
Let $G$ be an element of $K(X)$ with 
$v(G)=(H,f)-f$.
We set $v_D:=(H,C_0+D)-(C_0+D)$. 
Then $\deg_G(v)=0$ and $\chi_G(v)=-(1+(D,f))$.
We assume that $E_i$ are $\mu$-stable.
Let ${\frak g}'$ be the Lie algebra generated by
$C_1,\dots,C_n$.
By the remarks in section \ref{subsect:G-twisted},
we can construct an action of ${\frak g}'$
on $\bigoplus_D H_*(M_H^G(v_D),{\Bbb C})$, where
$D$ is an effective divisor with $(D,f)=0$.
\end{ex}



\subsection{Stable sheaves on an Enriques surface}

Let $X$ be an Enriques surface and $\pi:Y \to X$ be the covering $K3$ surface
of $X$.
Assume that $X$ contains a smooth $(-2)$ curve.
Let $C'$ be a connected component of $\pi^{-1}(C)$.
Let $H'$ be an ample divisor on $Y$ and set
$H:=\pi_*(H')$.
Then $H$ is an ample divisor on $X$ with
$(H,C)=2(H',C')$.
We take a semi-stable sheaf $E'$ on $Y$ with
$v(E')=(H',C')-C'$. $E$ is a rigid vector bundle.
If $H'$ is sufficiently ample, then 
Lemma \ref{lem:mu-stable} implies that
$E'$ is $\mu$-stable.
\begin{prop}
We set $E:=\pi_*(E')$. 
Then $E$ is a $\mu$-stable vector bundle
with the Mukai vector $(H,C)-C$ which satisfies $E \otimes K_X \cong E$
and
\begin{equation}
\begin{cases}
\Hom(E,E)={\Bbb C}\\
\Ext^1(E,E)=0\\
\Ext^2(E,E)={\Bbb C}.
\end{cases}
\end{equation}
\end{prop}
If there is a configuration of $(-2)$-curves, then 
as in the $K3$ surface case, we have an action of the Lie algebra 
associated to $(-2)$-curves on  
$\bigoplus_v H_*(M_H(v),{\Bbb C})$, where $v=(1+(D,H))+D+a \rho$,
$D \in \Pic(X)$, $a-1/2 \in {\Bbb Z}$.

\section{Actions associated to purely 1-dimensional exceptional sheaves}
\label{sect:1-dim}

\subsection{Purely 1-dimensional sheaves}
In this section, we shall consider Lie algebra actions
associated to purely 1-dimensional exceptional sheaves
such as line bundles on $(-2)$-curves.
Unfortunately we cannot construct the action
for the moduli spaces of stable torsion free
sheaves in general.
Instead, we can construct it for the moduli spaces
of purely 1-dimensional sheaves.
In some cases, the moduli spaces of stable torsion free sheaves
are deformation equivalent to moduli spaces of
purely 1-dimensional sheaves.
In this sense, we have an action
for the moduli spaces of stable torsion free
sheaves. This will be explained in \ref{subsect:rational}.
We also explain a partial result
on the moduli spaces of stable torsion free sheaves
in \ref{subsect:ADE}.

Let $(X,H)$ be a pair of a smooth projective surface $X$
and an ample divisor $H$ on $X$.

\begin{defn}\cite{Y:11}
Let $G$ be an element of $K(X)$ with $\rk G>0$.
A purely 1-dimensional sheaf $E$ is $G$-twisted stable, if 
\begin{equation}
\frac{\chi_G(F)}{(c_1(F),H)}<\frac{\chi_G(E)}{(c_1(E),H)}
\end{equation} 
for all proper subsheaf $F (\ne 0)$ of $E$.
\end{defn}
We have the following result whose proof is similar to
Lemma \ref{lem:key1}.

\begin{lem}\label{lem:key3}
Let $G$ be an element of $K(X)$ with $\rk G>0$ and $E_i$,
$i=1,2,\dots,s$, be purely 1-dimensional
$G$-twisted stable sheaves with $\chi_G(E_i)=0$.
Let $E$ be a purely 1-dimensional $G$-twisted stable sheaf with 
\begin{equation}\label{eq:min-cond3}
\chi_G(E)=\min\{\chi_G(E')>0| E' \in \Coh(X), \rk E'=0 \}
\end{equation}
or $E={\Bbb C}_P$, $P \in X$ with the condition \eqref{eq:min-cond3}.
\begin{enumerate}
\item[(1)]
Then every non-trivial extension
\begin{equation}
0 \to E_1 \to F \to E \to 0
\end{equation}
defines a $G$-twisted stable sheaf.
\item[(2)]
Let $V_i$ be a subspace of $\Hom(E_i,E)$.
Then  
$\phi:\bigoplus_{i=1}^s V_i \otimes E_i \to E$ is injective or surjective.
Moreover,
\begin{enumerate}
\item[(2-1)] 
if $\phi:\bigoplus_{i=1}^s V_i \otimes E_i \to E$ is injective, 
then the cokernel is a $G$-twisted stable purely 1-dimensional sheaf
or ${\Bbb C}_P$, $P \in X$,
\item[(2-2)]
if $\phi:\bigoplus_{i=1}^s V_i \otimes E_i \to E$ is surjective, then
$\ker \phi$ is $G$-twisted stable.
\end{enumerate}
\end{enumerate}
\end{lem}
\begin{lem}\label{lem:key4}
Let $G$ be an element of $K(X)$ with $\rk G>0$ and 
$E_i$, $i=1,2,\dots,s$, be purely 1-dimensional
$G$-twisted stable sheaves with $\chi_G(E_i)=0$.
Let $E$ be a purely 1-dimensional $G$-twisted stable sheaf with 
\begin{equation}\label{eq:min-cond4}
\chi_G(E)=\max\{\chi_G(E')<0| E' \in \Coh(X), \rk E'=0 \}.
\end{equation}
\begin{enumerate}
\item[(1)]
Then every non-trivial extension
\begin{equation}
0 \to E \to F \to E_1 \to 0
\end{equation}
defines a $G$-twisted stable sheaf.
\item[(2)]
Let $V_i$ be a subspace of $\Hom(E,E_i)$.
Then  
$\phi:E \to \bigoplus_{i=1}^s V_i^{\vee} \otimes E_i$ 
is injective or surjective.
Moreover,
\begin{enumerate}
\item[(2-1)] 
if $\phi:E \to \bigoplus_{i=1}^s V_i^{\vee} \otimes E_i$ is injective, 
then the
cokernel is a $G$-twisted stable purely 1-dimensional sheaf
or ${\Bbb C}_P$, $P \in X$,
\item[(2-2)]
if $\phi:E \to \bigoplus_{i=1}^s V_i^{\vee} \otimes E_i$ is surjective, 
then $\ker \phi$ is $G$-twisted stable.
\end{enumerate}
\end{enumerate}
\end{lem}

\begin{rem}\label{rem:H}
We set
\begin{equation}
d:=\min\{\chi_G(E')>0| E' \in \Coh(X), \rk E'=0 \}.
\end{equation}
For a purely 1-dimensional sheaf $E$ with
$\chi_G(E)=d$, $E$ is $G$-twisted stable if and only if
$\chi_G(F) \leq 0$ for all proper subsheaf $F$ of $E$.
Thus the $G$-twisted stability does not depend on the choice of 
$H$.
\end{rem}

\begin{defn}\label{defn:v}
For a complex ${\Bbb E}$ with $\rk({\Bbb E})=0$,
we set
\begin{equation}
v({\Bbb E}):=(c_1({\Bbb E}),\chi({\Bbb E})) 
\in H^2(X,{\Bbb Z}) \times {\Bbb Z}.
\end{equation}
We define a pairing of 
$v_i:=(\xi_i,a_i) \in H^2(X,{\Bbb Z}) \times {\Bbb Z}$, $i=1,2$ by
\begin{equation}
 \langle v_1,v_2 \rangle:=(\xi_1,\xi_2) \in {\Bbb Z}.
\end{equation}
\end{defn}
Then the Riemann-Roch theorem says that
\begin{equation}
\chi({\Bbb E},{\Bbb F})=-\langle v({\Bbb E}),v({\Bbb F}) \rangle
\end{equation}
for ${\Bbb E},{\Bbb F} \in {\bf D}(X)$ with
$\rk({\Bbb E})=\rk ({\Bbb F})=0$. 
We set $\rho:=v({\Bbb C}_P)=(0,1)$.

\begin{defn}
Let ${\Bbb E} \in {\bf D}(X)$ be an object such that 
$\rk ({\Bbb E})=0$ and
\begin{equation}
\chi_G({\Bbb E})=\min\{\chi_G({\Bbb E}')>0 |{\Bbb E}' \in {\bf D}(X),
\rk({\Bbb E}')=0 \}.
\end{equation}
${\Bbb E}$ is $G$-twisted stable, if
\begin{enumerate}
\item
$H^i({\Bbb E})=0$, $i \ne 0$ and $H^0({\Bbb E})$ is $G$-twisted stable, or 
\item
$H^i({\Bbb E})=0$, $i \ne -1$ and $H^{-1}({\Bbb E})$ 
is $G$-twisted stable.
\end{enumerate}
Let $M_H^G(v)$ be the moduli space of $G$-twisted stable complexes
${\Bbb E}$ with $v({\Bbb E})=v$.
\end{defn}

Let $E_i$, $i=1,2,\dots,n$ be $G$-twisted stable purely 1-dimensional
sheaves such that $\chi_G(E_i)=0$,
$E_i \otimes K_X \cong E_i$ and
$\langle v(E_i)^2 \rangle=-2$.
We set $v_i:=v(E_i)$.
Let ${\frak g}$ be the Lie algebra associated to $E_i$, $i=1,\dots,n$.
By using Lemma \ref{lem:key3}, \ref{lem:key4}, we get the following 
similar results to the results in section \ref{sect:action}.
\begin{prop}\label{prop:action-1dim}
For all ${\Bbb E} \in M_H^G(v+\sum_i x_i v_i)$, $x_i \in {\Bbb Z}$,
we assume that 
\begin{equation}
\chi_G({\Bbb E})=\min\{\chi_G({\Bbb E}')>0| {\Bbb E}' \in {\bf D}(X), 
\rk {\Bbb E}'=0 \}
\end{equation}
and ${\Bbb E}$ satisfies \eqref{eq:smooth-condition}.
Then we have an action of ${\frak g}$ on
$\bigoplus_{x_i \in {\Bbb Z}} H_*(M_H^G(v+\sum_i x_i v_i),{\Bbb C})$.
\end{prop}

Let $C$ be an irreducible $(-2)$-curve on $X$.
If $G={\cal O}_X$, then ${\cal O}_C(-1)$ is a stable sheaf
with $\chi({\cal O}_C(-1))=0$.
Then we can apply Proposition \ref{prop:action-1dim}. 
\begin{cor}
Let $X$ be a K3 surface. Then $M_H((D,1)) \ne \emptyset$
for all $H$.
\end{cor}

\begin{proof}
By Proposition \ref{prop:action-1dim}, we have isomorphims
\begin{equation}
H_*(M_H((D,1)),{\Bbb C}) \cong
H_*(M_H((D+(D,C)C,1)),{\Bbb C})
\end{equation}
for all irreducible
$(-2)$-curves $C$.
Hence we can reduce the proof to the case where
$D$ is a nef divisor or $D$ is a smooth rational curve.
If $D$ is nef, then \cite[Rem. 3.4]{Y:12} implies that
$M_H((D,1)) \ne \emptyset$.
Therefore our claim holds.
\end{proof}
 
\begin{rem}
We can show that $M_H^G((D,n)) \ne \emptyset$
for a general $(H,G)$ by a different method.     
\end{rem}

\begin{lem}\label{lem:half-K3}
Let $X$ be a 9 point blow-up of ${\Bbb P}^2$
and assume that $|-K_X|$ contains a reducible curve
$Y=\sum_{i=0}^n a_i C_i$, where $C_i$ are smooth $(-2)$-curves.
Then every $G$-twisted stable purely 1-dimensional sheaf $E$
with $(c_1(E),K_X)<0$
satsifies \eqref{eq:smooth-condition}.
\end{lem}

\begin{proof}
Assume that there is a non-zero map
$\psi:E \to E(K_X)=E(-Y)$.
By the homomorphism ${\cal O}_X(-Y) \to {\cal O}_X$,
we have a homomorphism 
$E(-Y) \to E$, which is isomorphic on
$\Div(E) \setminus Y \neq \emptyset$.
If $E \to E(-Y) \to E$ is a zero map,
then
$F:=\psi(E)(-K_X)$ satisfies
$\Supp(F(K_X)) \subset Y$ and
$\chi_G(E)/(c_1(E),H)< \chi_G(F(K_X))/(c_1(F),H)$.
On the other hand, since $F$ is a proper subsheaf of $E$,
we have $\chi_G(F)/(c_1(F),H)<\chi_G(E)/(c_1(E),H)$.
Since $(C_i,K_X)=0$, we get $(c_1(F),K_X)=0$.
This means that $\chi_G(F(K_X))=\chi_G(F)$.
Then we get $\chi_G(E)/(c_1(E),H)<\chi_G(E)/(c_1(E),H)$.
This is a contradiction.
Therefore $E \to E(-Y) \to E$ is a non-zero map.
Then by using the stability of $E$ and $(\Div(E),Y)>0$, 
we get a contradiction.
Hence we conclude that $\Hom(E,E(K_X))=0$.
\end{proof}   

\begin{cor}
Under the assumption in Lemma \ref{lem:half-K3},
we have an action of the affine Lie algebra
associated to $C_i$, $0 \leq i \leq n$
on $\bigoplus_D H_*(M_H((D,1)),{\Bbb C})$, where $(D,K_X)<0$.
\end{cor}

\begin{prop}\label{prop:(-2)}
Let $C_i$, $i=0,1,\dots, n$ 
be a configulation of smooth $(-2)$-curves of $ADE$ or affine type such that
$K_X$ is trivial in a neighbouhood of $\cup_i C_i$ and
$(C_i,C_j) \leq 1$ for $i \ne j$. 
Let $D:=\sum_{i=0}^n b_i C_i, b_i>0$ be an effective divisor
such that $(D^2)=-2$.
There is a $G$-twisted stable sheaf $E$ with
$(c_1(E),\chi(E))=(D,m)$ for a general $(H,G)$.
\end{prop}

\begin{proof}
If $n=0$, then $D=C_0$ and obviously the claim holds.
Hence we may assume $n>0$ and $\cup_i C_i$ is connected.
We set $v_i:=v({\cal O}_{C_i}(-1))$.
We first show that $M_H(\rho+\sum_i b_i v_i) \ne \emptyset$.
Assume that there is a stable sheaf
$E$ such that $\Supp(E) \subset \cup_i C_i$.
Since $K_X$ is trivial in a neighbourhood of $\cup_i C_i$, 
$\Ext^2(E,E) \cong \Hom(E,E \otimes K_X)^{\vee}={\Bbb C}$.
Hence we see that $(c_1(E)^2) \geq -2$ and the equality holds 
when $\Ext^1(E,E)=0$.
In particular $M_H(\rho+\sum_i b_i v_i)$ is smooth.
Let $R_{v_i}$ be the $(-2)$-reflection defined by $v_i$ and 
$W$ the Weyl group generated by $R_{v_i}$, $i=0,1,\dots,n$.
Then by the action of $W$,
we have an isomorphism
$M_H(\rho+\sum_i b_i v_i) \to M_H(\rho+v_j)$ for some $j$.
Obviously $M_H(\rho+v_j)=\{{\cal O}_{C_j} \}$.
Therefore $M_H(\rho+\sum_i b_i v_i) \ne \emptyset$.

We shall treat the general cases.
Since $K_X$ is trivial in a neighbouhood of $\cup_i C_i$, 
by using \cite[Prop. 2.7]{Y:11},
we see that the non-emptyness 
of $M_H^G(v)$ does not depend on the choice of a general
$(H,G)$.
There is a divisor $C$ such that $(C,D)=1$.
Indeed we take an element $w \in W$ 
such that $w(D)=C_i$ for some $i$.
Then $1=(C_j,C_i)=(w(C_j),D)$ for a $j$.  
Let $E$ be a stable sheaf with
$(c_1(E),\chi(E))=(D,1)$. 
Then $E(nC)$ is a ${\cal O}_X(nC)$-twisted stable sheaf
with $\chi(E(nC))=1+n$.
Therefore our claim holds for general cases.

\end{proof}

\begin{ex}\label{ex:RDP}
Let $Y$ be a germ of a rational double point and 
$\pi:X \to Y$ the minimal resolution.
Let $H$ be a $\pi$-ample divisor on $X$.
Let $C_i$, $i=1,2,\dots,n$ be irreducible components
of the exceptional divisor.
We set $v_i:=v({\cal O}_{C_i}(-1))$.
Let ${\frak g}$ be the Lie algebra associated to $C_i$.
We note that $K_X \cong {\cal O}_X$.
For a coherent sheaf $E$ on $X$ with a compact support,
we can define the stability with resppect to $H$. 
For a stable sheaf $E$ with
$v(E)=\rho+\sum_i n_i v_i$, $\dim \Ext^1(E,E)=
\langle (\rho+\sum_i n_i v_i)^2 \rangle+2$.
Hence we get 
\begin{equation}
\dim \Ext^1(E,E)=
\begin{cases}
2, \text{ if $v=\rho$},\\
0, \text{ if $\langle (\sum_i n_i v_i)^2 \rangle=-2$}.
\end{cases}
\end{equation}
If $v=\rho$, then all stable sheaves are of the form
${\Bbb C}_P, P \in X$.
Hence $M_H(\rho)$ has a coarse moduli space which is isomorphic to
$X$.
Hence $M_H(\rho)$ is smooth.
If $\langle (\sum_i n_i v_i)^2 \rangle=-2$, then 
the proof of Proposition \ref{prop:(-2)} implies that
 $M_H(\rho+\sum_i n_i v_i)$
is not empty and consists of a stable sheaf on the exceptional divisors.     
Then we have an action of ${\frak g}$ on
$\bigoplus_{n_i \in {\Bbb Z}} H_*(M_H(\rho+\sum_i n_i v_i),{\Bbb C})$.
Indeed the submodule consisting of the middle degree homology groups
is isomorphic to ${\frak g}$.
For the structure of $M_H(\rho+\sum_i n_i v_i)$, we get the following:
Let $D=\sum_i n_i C_i$ be an effective divisor with $(D^2)=-2$.
Then $M_H(\rho+\sum_i n_i v_i)=\{{\cal O}_D\}$
and $M_H(\rho-\sum_i n_i v_i)=\{{\cal O}_D(D)\}$.

Proof of the claim:
We note that $\chi({\cal O}_D)=-(D^2)/2=1$.
If there is a quotient ${\cal O}_D \to {\cal O}_{D'}$,
then since $({D'}^2)<0$, we have 
$\chi({\cal O}_{D'})=-({D'}^2)/2 \geq 1$.  
Therefore ${\cal O}_D$ is stable.
We note that ${\cal O}_D(D)$ is the derived dual of
${\cal O}_D$.
By using this fact, we can easily see the stability of 
${\cal O}_D(D)$. 
\end{ex}

\begin{ex}
Let $C$ be a germ of a curve at $P$
and 
$\pi:X \to C$ an elliptic surface with a section
$\sigma$.
Let $H$ be a $\pi$ ample diviosr on $X$.
Assume that $\pi^{-1}(P)$ is reducible and consists of 
smooth $(-2)$-curves $C_i$, $i=0,1,\dots,n$:
$\pi^{-1}(P)=\sum_{i=0}^n a_i C_i$.
We may assume that $a_0=1$ and $(\sigma,C_0)=1$.
We set $v_i:=v({\cal O}_{C_i}(-1))$.
Then we see that $M_H(\rho+\sum_i n_i v_i)$ is smooth with
\begin{equation}
\dim M_H(\rho+\sum_i n_i v_i)=\langle (\sum_i n_i v_i)^2 \rangle+2
\end{equation}
and $M_H(\rho+\sum_i n_i v_i) \cong X$, 
if $\langle (\sum_i n_i v_i)^2 \rangle=0$.
We also have an action of affine Lie algebra ${\frak g}$
associated to
$v_i$ on $\bigoplus_{n_i \in {\Bbb Z}}
H_*(M_H(\rho+\sum_i n_i v_i),{\Bbb C})$. 
Indeed, if $(C_i,C_j) \leq 1$ for all $i \ne j$, 
then the result obviously holds.
If $(C_i,C_j)=2$, then we can directly check the commutation relation
\eqref{eq:relation3}.
We set $\delta:=v({\cal O}_{\pi^{-1}(P)})
=\sum_{i=0}^n a_i v_i$.  
If $ \sum_i n_i v_i=m \delta$, $m \in {\Bbb Z}$,
then under an identification $M_H(\rho+m\delta) \cong X$,
we have an isomorphism
\begin{equation}
H_2(M_H(\rho+m\delta),{\Bbb C}) \cong 
{\Bbb C}[\sigma] \oplus \bigoplus_{i=1}^n {\Bbb C}[C_i].
\end{equation}
%
Let 
\begin{equation}
{\frak g}={\Bbb C}[t,t^{-1}] \otimes \overline{\frak g} 
\oplus {\Bbb C}c \oplus {\Bbb C}d
\end{equation}
be the standard expression of the affine Lie algebra,
where $c$ is the center of ${\frak g}$.
Then we have an exact sequence of ${\frak g}$-modules:
\begin{equation}
0 \to
{\Bbb C}[t,t^{-1}] \otimes \overline{\frak g}
\to
\bigoplus_{n_i \in {\Bbb Z}}
H_{mid}(M_H(\rho+\sum_i n_i v_i),{\Bbb C}) 
\to {\Bbb C}[t,t^{-1}] \to 0.
\end{equation}
where $H_{mid}(\star)$ is the middle degree homology group of
$\star$ and ${\Bbb C}[t,t^{-1}]$ is the
${\frak g}/[{\frak g},{\frak g}]={\Bbb C}d$-module.
%
\end{ex}

\subsection{Moduli of stable sheaves on elliptic
surfaces}\label{subsect:elliptic}

We first collect some basic facts on the moduli spaces
of stable sheaves on elliptic surfaces $X$.
If $K_X$ is not numerically trivial, then
we do not have a good invariant of
torsion free sheaf $E$ which is a suitable generalization
of the Mukai vector.
In these cases,  we shall use 
$\gamma(E):=(\rk(E),c_1(E),\chi(E)) \in H^*(X,{\Bbb Z})$ as an invariant
of $E$.
If $\rk E=0$, then $\gamma(E)=(0,c_1(E),\chi(E))$
is the same as the Mukai vector $v(E)$
defined in Definition \ref{defn:v}.
We denote the moduli space
of $G$-twisted stable sheaves $E$ 
on $X$ with $\gamma(E)=(r,\xi,\chi)$
by $M_H^G(r,\xi,\chi)$. We also denote the moduli of
$G$-twisted semi-stable sheaves by $\overline{M}_H^G(r,\xi,\chi)$.

Let $\pi:X \to C$ be an elliptic surface.
Let $f$ be a smooth fiber and
$L$ a nef and big divisor on $X$.
Since $(L,f)>0$, replacing $L$ by $L+nf$, $n>0$,
we may assume that $(L,C')>0$ unless 
$C'$ is a $(-2)$-curve in a fiber of $\pi$. 
Let $G$ be a locally free sheaf on $X$ such that
$\rk G=r$ and $c_1(G)=d\sigma$ with $\gcd(r,d)=1$.
We first study the stability condition, 
when the polarization is sufficiently close to $f$.

\begin{lem}\label{lem:elliptic}
For $(\xi,\chi) \in \NS(X) \times {\Bbb Z}$ with $r\chi-(c_1(G),\xi)>0$, 
we take $(n, \varepsilon) \in {\Bbb Z} \times \NS(X) \otimes {\Bbb Q}$
such that $n \gg 0$ and 
$\varepsilon$ is an ample ${\Bbb Q}$-divisor with 
$|\varepsilon| \ll 1$.
Let $E$ be a purely 1-dimensional sheaf with 
$v(E)=(\xi,\chi)$. 
\begin{enumerate}
\item
$E$ is $G$-twisted stable with respect to 
$L+nf+\varepsilon$ if and only if
for all proper subsheaf $F$ of $E$, one of the following holds
\begin{enumerate}
\item
$\chi_G(E)/(c_1(E),f)>\chi_G(F)/(c_1(F),f)$ or
\item
$\chi_G(E)/(c_1(E),f)=\chi_G(F)/(c_1(F),f)$ and 
$\chi_G(E)/(c_1(E),L)>\chi_G(F)/(c_1(F),L)$ or
\item
$\chi_G(E)/(c_1(E),f)=\chi_G(F)/(c_1(F),f)$,
$\chi_G(E)/(c_1(E),L)=\chi_G(F)/(c_1(F),L)$ and 
$\chi_G(E)/(c_1(E),\varepsilon)>\chi_G(F)/(c_1(F),\varepsilon)$.
\end{enumerate}
\item
Moreover if $L+nf$ is ample and $\gcd((c_1(E),f),(c_1(E),L),\chi_G(E))=1$, 
then 
there is no properly $G$-twisted semi-stable sheaf $E$ with respect to
$L+nf$.
In particular, $E$ is $G$-twisted stable with respect to $L+nf$ 
if and only if $E$ is $G(\eta)$-twisted stable with respect to 
$L+nf+\varepsilon'$,
where $\eta, \varepsilon'$ are sufficiently small ${\Bbb Q}$-divisors.  
\end{enumerate}
\end{lem}

\begin{proof}
We note that $\chi_G(E)=r\chi-(c_1(G),\xi)>0$.
If $\chi_G(E)(c_1(F),f)-\chi_G(F)(c_1(E),f)<0$, 
then we see that
\begin{equation}
\begin{split}
& \chi_G(E)(c_1(F),L+nf+\varepsilon)-\chi_G(F)(c_1(E),L+nf+\varepsilon)\\
\leq & n(\chi_G(E)(c_1(F),f)-\chi_G(F)(c_1(E),f))+
\chi_G(E)(c_1(E),L+\varepsilon).
\end{split}
\end{equation}
If $\chi_G(E)(c_1(F),f)-\chi_G(F)(c_1(E),f) \geq 0$ and
$\chi_G(F) \geq 0$, then 
\begin{equation}
\begin{split}
& \chi_G(E)(c_1(F),L+nf+\varepsilon)-\chi_G(F)(c_1(E),L+nf+\varepsilon)\\
\geq & n(\chi_G(E)(c_1(F),f)-\chi_G(F)(c_1(E),f))\\
& \quad +
\frac{\chi_G(E)}{(c_1(E),f)}
\left((c_1(E),f)(c_1(F),L+\varepsilon)-(c_1(F),f)(c_1(E),L+\varepsilon)
\right)\\
\geq & n(\chi_G(E)(c_1(F),f)-\chi_G(F)(c_1(E),f))-
\chi_G(E)(c_1(E),L+\varepsilon).
\end{split}
\end{equation}
By using these inequalities, we can show claim (i). 
Moreover, if $L+nf$ is ample,
then the equalities
\begin{equation}
\begin{split}
\chi_G(E)(c_1(F),f)-\chi_G(F)(c_1(E),f)& =0,\\
\chi_G(E)(c_1(F),L)-\chi_G(F)(c_1(E),L)& =0
\end{split}
\end{equation}
imply that
$(c_1(F),f)=(c_1(F),L)=\chi_G(F)=0$ or
$(c_1(E/F),f)=(c_1(E/F),L)=\chi_G(E/F)=0$.
By the ampleness of $L+nf$, we get
$F=0$ or $E/F=0$. 
Therefore the claim (ii) holds. 
\end{proof}

\begin{lem}
Let $\pi:X \to C$ be an elliptic surface and
$f$ a fiber of $\pi$.
If $(D,f)=1$, then ${\Bbb E} \in M_H(0,D,1)$ satisfies
\eqref{eq:smooth-condition}.
\end{lem}
For the proof, see \cite[Prop. 3.18]{Y:11}.
%
Let $\pi^{-1}(p)=\sum_{i=0}^n a_i C_i$ be a singular fiber of $\pi$ such that
$(C_i,C_j) \leq 1$.
We may assume that $(C_0,\sigma)=a_0=1$.
\begin{lem}\label{lem:E0}
There is a $G$-twisted stable sheaf $E_0$ with
$c_1(E_0)=(r-1)f+C_0$ and $\chi_G(E_0)=0$.
\end{lem}

\begin{proof}
By Proposition \ref{prop:(-2)}, there is a
$G$-twisted semi-stable sheaf $E_0$ with
$c_1(E_0)=(r-1)f+C_0$ and $\chi_G(E_0)=0$.
Assume that $E_0$ is $S$-equivalent to
$\bigoplus_i F_i$, where $F_i$ are $G$-twisted stable sheaves with
$\chi_G(F_i)=0$.
Since $\Supp(F_i)$ does not contain $\sigma$,
$(c_1(F_i),\sigma) \geq 0$.
Since $\chi_G(F_i)=r\chi(F_i)-d(\sigma,c_1(F_i))$,
there is an integer $i_0$ such that
$(\sigma,c_1(F_{i_0}))=r$ and $(\sigma,c_1(F_{i}))=0$ for
$i \ne i_0$.
Thus $\Supp(F_i)$, $i \ne i_0$ do not contain $C_0$, which implies that
$(c_1(F_{i}),C_0) \geq 0$, $i \ne i_0$.
Then we see that 
$(c_1(F_{i_0})^2) \leq (c_1(E_0)^2)+
((\sum_{i \ne i_0} c_1(F_i))^2)<-2$.
This is a contradiction.
Therefore $E_0$ is $G$-twisted stable.
\end{proof} 

\begin{lem}
We set $E_i:={\cal O}_{C_i}(-1)$, $i>0$.
Let $E$ be a properly $G$-twisted semi-stable sheaf with
$\gamma(E)=(0,rf,d)$ and $\Supp(E)=\sum_i a_i C_i$.
Then $E$ is $S$-equivalent to $\bigoplus_i E_i^{\oplus a_i}$.
\end{lem}

\begin{proof}
By the proof of Lemma \ref{lem:E0},
it is sufficient to prove the following.
Let $F$ be a $G$-twisted stable sheaf with
$c_1(F)=\sum_{i>0} n_i C_i$ and
$\chi(F)=\chi_G(F)=0$. Then $F=E_i$ for an $i>0$:

Since $(c_1(F)^2)<0$,
we can choose an integer $i$ such that $(c_1(F),C_i)<0$.
Then $\chi(E_i,F)>0$, which implies that
$\Hom(E_i,F) \ne 0$ or $\Hom(F,E_i) \ne 0$.
By the stability of $E_i$ and $F$, 
we see that $E_i \cong F$.
Thus the claim holds. 
\end{proof}

We take a sufficiently small ${\Bbb Q}$-divisor $\eta$ such that
$(\sigma,\eta)=(f,\eta)=0$ and $\chi_{G(\eta)}(E_i)<0$ for $i>0$.
Then 
in the same way as in \cite{O-Y:1}, we see that 
$Y:=M_H^{G(\eta)}(0,rf,d)$ is a resolution of $\overline{M}_H^G(0,rf,d)$
at $\bigoplus_i E_i^{\oplus a_i}$ and the exceptional divisors
are 
\begin{equation}
C_i':=\{E \in Y|\Ext^2(E,E_i) \ne 0 \} \cong {\Bbb P}^1, i>0
\end{equation}
and $(C_i',C_j')=(C_i,C_j)$.

Let ${\frak g}$ (resp. $\overline{\frak g}$)
be the affine Lie algebra associated to $E_i$,
$i \geq 0$ (resp. the finite Lie algebra associated to $E_i$,
$i \geq 1$). 

\begin{prop}\label{prop:elliptic}
Let $\pi:X \to C$ be an elliptic surface with a section
$\sigma$ as above.
Assume that (1) $l=1$, or (2) $X$ is rational or of type $K3$. 
We set $L:=\sigma+(1-(\sigma^2))f$.
Let $G$ be a locally free sheaf on $X$ such that
$\rk G=r$ and $c_1(G)=d\sigma$ with $\gcd(r,d)=1$.
Then $\overline{\frak g}$ acts on 
$\bigoplus_D H_*(M_{L+nf+\varepsilon}^G(0,l\sigma+D,k),{\Bbb C})$,
where $D$ is an effective divisor on fibers
with $(l\sigma+D)^2+p_g+1 \geq 0$ and $k$ is an integer
with 
$\chi_G:=rk-(c_1(G),l\sigma+D)>0$ and 
$\gcd(l,(D,\sigma),\chi_G)=1$.
Moreover we also have an action of ${\frak g}$, if
$\gcd(l,\chi_G)=1$. 
\end{prop}

\begin{proof}
We first note that \eqref{eq:smooth-condition} holds,
under (1) or (2).
We first note that $(c_1(E)\pm c_1(E_i),\sigma)=(c_1(E),\sigma)$
for $i>0$.
If $\gcd(\chi_G(E),(c_1(E),f),(c_1(E),\sigma))=1$, then
the statements in Lemma \ref{lem:key3} hold,
where $E_1$ in Lemma \ref{lem:key3} corresponds to $E_i, i>0$. 
Hence we get our claim for $\overline{\frak g}$.
Moreover if $\gcd(\chi_G(E),(c_1(E),f))=1$, then we can apply the results
in Lemma \ref{lem:key3} for $E_0$.
Therefore our claim also holds for ${\frak g}$.
\end{proof}

\begin{cor}\label{cor:weyl}
Under the same notations as above,
the Poincar\'{e} polynomial $P(M_{L+nf+\varepsilon}^G(0,l\sigma+D,k))$
is $W({\frak g})$-invariant:
\begin{equation}
P(M_{L+nf+\varepsilon}^G(0,w(l\sigma+D),k))=
P(M_{L+nf+\varepsilon}^G(0,l\sigma+D,k)),\;w \in W({\frak g}),
\end{equation}
where $W({\frak g})$ is the Weyl group of ${\frak g}$.
\end{cor}
Let $X$ be a rational elliptic surface with a singular fiber of type 
$E_8^{(1)}$. 
As we shall see in subsection \ref{subsect:rational},
$M_{L+nf+\varepsilon}^G(0,l\sigma+D,k)$ is related to a moduli 
space of torsion free sheaves. 
In \cite{MNWV}, \cite{Y:4} and \cite{Iq},
it is observed that the Eular characteristic of
$M_{L+nf+\varepsilon}^G(0,l\sigma+D,k)$ is Weyl group invariant.
Proposition \ref{prop:elliptic} gives an explanation of
this invariance.

\subsection{Moduli of stable sheaves on rational elliptic
surfaces}\label{subsect:rational}

Let $\pi:X \to {\Bbb P}^1$ be a rational elliptic surface with
a section $\sigma$. 
Then there is a family of elliptic surfaces
$\widetilde{\pi}:{\cal X} \to {\Bbb P}^1_T$ over a scheme $T$ such that
\begin{enumerate}
\item
${\cal X}_{t_0} \cong X$, $t_0 \in T$,
\item
there is a section $\widetilde{\sigma}$ of $\widetilde{\pi}$
with $\widetilde{\sigma}_{t_0}=\sigma$
and 
\item
for a general point $t \in T$,
${\cal X}_t$ is a nodal elliptic surface, that is,
all singular fibers are irreducible nodal curves.
\end{enumerate}
Let $T_0$ be the open subset of $T$ consisting of nodal
elliptic surfaces.
Replacing $T$ by a suitable covering of $T$, we may assume that
$\Pic({\cal X}/T) \cong R^2\phi_*({\Bbb Z})$ is a trivial local system,
where $\phi:{\cal X} \to T$ is the projection.
\begin{NB}
Let $\xi_i$, $i=1,2,\dots,10$ be a basis of $H^2(X,{\Bbb Z})$.
For each $\xi_i \in H^2(X,{\Bbb Z})$,
let $\Pic^{\xi_i}({\cal X}/T) \to T$ be the relative moduli space
of line bundles which is a deformation of a line bundle
$L_i$ on $X$ with $c_1(L_i)=\xi_i$.
It is \'{e}tale over $T$. We have a projective bundle
$P_i \to \Pic^{\xi_i}({\cal X}/T)$ and a family of line bundles
${\cal L}_i$ on $P_i \times_T {\cal X}$, which gives a section
$c_1({\cal L}_i)$ of $R^2\phi_*({\Bbb Z})$ over $P_i$.
Then we have a trivialization of $R^2\phi_*({\Bbb Z})$ on
$P_1 \times_T P_2 \times \dots \times_T P_{10}$.
\end{NB}   
Hence there is a relatively ample ${\Bbb Q}$-divisor ${\cal H}$
on ${\cal X}$.
Moreover, by adding $m\sigma+nf$, we may assume that
${\cal H}=nf+m\sigma+\varepsilon$,
$\varepsilon \in ({\Bbb Q}\sigma+{\Bbb Q}f)^{\perp}$,
$n \gg m \gg |(\varepsilon^2)|$, where
we use the identification 
$R^2\phi_*{\Bbb Z} \cong H^2(X,{\Bbb Z})$.
For positive integers $r,d$ with $\gcd(r,d)=1$,
we take a vector bundle $G$ of rank $r$ and $c_1(G)=d\sigma$
on ${\cal X}$.
We set $\gamma:=(r,\xi,\chi) \in R^* \phi_* {\Bbb Z}$.
Then we have a family of moduli spaces of semi-stable sheaves
$\psi:\overline{M}_{({\cal X},{\cal H})/T}^G(\gamma) \to T$, which is smooth
on the locus of stable sheaves.
Assume that $\gamma$ is primitive and $G$ is general with respect to
$\gamma$.
Then $\overline{M}_{({\cal X},{\cal H})/T}^G(\gamma)$ 
consists of stable sheaves.
%
%
For the existence of stable sheaves, see Appendix \ref{subsect:exist}. 

From now on, we assume that $\gamma=(0,\xi,\chi)$, where
$\xi=l\sigma+kf+D$, $l>0$, $\gcd(l,(\xi,\sigma),r\chi-(c_1(G),\xi))=1$ 
and $(D,f)=(D,\sigma)=0$. 
We take a sufficiently small ${\Bbb Q}$-divisor $\eta$
such that $\chi_{G(\eta)}({\cal O}_{C_i}(-1))<0$ for all $i>0$.
We set ${\cal Y}:=M_{({\cal X},{\cal H})/T}^{G(\eta)}(0,r \sigma, d)$.
Then ${\cal Y}_t$, $t \in T$ are smooth projective surface isomorphic to
${\cal X}_t$. Hence ${\cal Y} \to T$ is a smooth morphism.
We have an isomorphism ${\cal Y} \times_T T_0 \cong {\cal X} \times_T T_0$
over $T_0$ (cf. \cite{Y:11}).
Let ${\cal H}'$ be a relatively ample ${\Bbb Q}$-divisor on ${\cal Y}$
whose restriction to ${\cal Y} \times_T T_0$ corresponds to
a divisor on ${\cal X} \times_T T_0$ which is sufficiently close to
$m\sigma+nf$.
By Lemma \ref{lem:elliptic}, we have an isomorphism
\begin{equation}
M_{({\cal X} \times_{T}T_0,{\cal H}_{|\phi^{-1}(T_0)})/T_0}^G(\gamma) 
\cong M_{({\cal X}\times_{T}T_0,m\sigma+nf)/T_0}^{G(-\varepsilon)}(\gamma).
\end{equation}
By our assumption, there is a universal family
${\cal E}$ on ${\cal X} \times_T {\cal Y}$.
We consider a family of Fourier-Mukai transforms
$\Phi_{{\cal X} \to {\cal Y}}^{{\cal E}^{\vee}}:
{\bf D}({\cal X}) \to {\bf D}({\cal Y})$.
If $r\chi-(c_1(G),\xi)>0$, then
$\Phi_{{\cal X} \to {\cal Y}}^{{\cal E}^{\vee}}$ induces a
birational map
\begin{equation}
\zeta:M_{({\cal X},{\cal H})/T}^G(\gamma)
\cdots \to M_{({\cal Y},{\cal H}')/T}^{{G'}^{\vee}}(\gamma')
\end{equation}
which is an isomorphism over $T_0$,
where
$G':=\Phi_{{\cal X} \to {\cal Y}}^{{\cal E}^{\vee}}({\cal O}_{\sigma})$
(see \cite[Thm. 3.13, Rem. 3.1]{Y:11}).
Let ${\cal Z}$ be the graph of this birational correspondence.
Then the cycle $[{\cal Z}]_{t_0}$ induces an isomorphism
of the homology groups 
\begin{equation}
H_*(M_{{\cal H}_{t_0}}^G(\gamma),{\Bbb Z}) \to 
H_*(M_{{\cal H}_{t_0}'}^{{G'}^{\vee}}(\gamma'),{\Bbb Z})
\end{equation}
via the convolution product. 
Let $E_i$, $i=0,1,\dots,n$ 
be $G$-twisted stable sheaves on $X$ in section \ref{subsect:elliptic}.
We set $Y:={\cal Y}_{t_0}$ and 
\begin{equation}
\begin{split}
\rho:=&\Phi_{{\cal X}_{t_0} \to {\cal Y}_{t_0}}^{{\cal E}_{t_0}^{\vee}}
({\Bbb C}_x)) \in 
K(Y), x \in X\\
u_i:=&\Phi_{{\cal X}_{t_0} \to 
{\cal Y}_{t_0}}^{{\cal E}_{t_0}^{\vee}}(E_i)) \in 
K(Y), i=0,1,\dots,n.
\end{split}
\end{equation}
Then $\sum_{i=0}^n a_i u_i=
\Phi_{{\cal X}_{t_0} \to 
{\cal Y}_{t_0}}^{{\cal E}_{t_0}^{\vee}}(\bigoplus_{i=0}^n E_i^{\oplus a_i})
={\Bbb C}_y, y \in Y$.
By Proposition \ref{prop:elliptic}, we get the following:
\begin{prop}
We have an action of $\overline{\frak g}$ on the homology groups 
\begin{equation}
\bigoplus_{n_i,k} 
H_*(M_H^{{G'}^{\vee}}(\gamma(lG'+\sum_i n_i u_i+k\rho)),
n_i \in {\Bbb Q},
\end{equation}
where $\sum_i n_i u_i \in K(Y)$,
$H$ is sufficiently close to $f$,
$l\chi_G({\cal O}_{\sigma})+kr>0$ and 
$\gcd(l,rn_0,kr)=1$.
Moreover if $\gcd(l,kr)=1$, 
then we have an action of ${\frak g}$.
\end{prop}

\begin{proof}
We note that $(c_1(\sum_i n_i u_i),\sigma)=n_0 r$
and $\chi_G(l{\cal O}_{\sigma}+k{\Bbb C}_x) \equiv kr \mod l$.
Hence the claim holds.
\end{proof}

%
%


\subsection{Moduli of stable vector bundles on an $ADE$-configuration.}
\label{subsect:ADE}

In this subsection, we explain a relation 
with a paper by Nakajima \cite{Na:gauge}. 
Let $X$ be a smooth projective surface containing
a $ADE$-type configuration of smooth rational curves $C_i$, $i=1,2,\dots,n$.
Assume that there is a nef and big divisor $H$ such that
$(C_i,H)=0$ for all $C_i$.

For $\xi \in \NS(X)$ and $d \geq 0$,
we set
\begin{equation}
B_{(\xi,d)}:=\left\{ \left. x \in \oplus_{i=1}^n {\Bbb Z}C_i \right| 
(x^2)-2(\xi,x)+d \geq 0 \right\}.
\end{equation}
Since $\bigoplus_{i=1}^n {\Bbb Z}C_i$ is negative definite,
$B_{(\xi,d)}$ is a finite set.
Let $r$ be a positive integer such that
$2r>(x^2)-2(\xi,x)+d$ for all $x \in B_{(\xi,d)}$.
Assume that there is an integer $\chi_0$ 
such that $d=(\xi^2)-2r \chi_0-r(K_X,\xi)+(r^2+1)\chi({\cal O}_X)$.
\begin{defn}
Let $M_H(r,\xi+y,\chi)^{\mu}$, $y \in \bigoplus_{i=1}^n {\Bbb Z}C_i,
\chi \in {\Bbb Z}$ 
be the moduli space of $\mu$-stable sheaves
$E$ with respect to $H$ such that $\gamma(E)=(r,\xi+y,\chi)$.
\end{defn}
$M_H(r,\xi+y,\chi)^{\mu}$ is contained in a moduli space
of $\mu$-stable sheaves with respect to an ample
divisor $H'$ which is sufficiently close to $H$.
If $\gcd(r,(\xi,H))=1$, then 
$M_H(r,\xi+y,\chi)^{\mu}$ is projective.
We assume that 
\eqref{eq:smooth-condition} holds for all 
$E \in M_H(r,\xi+y,\chi)^{\mu}$, $y \in \bigoplus_{i=1}^n {\Bbb Z}C_i,
\chi \in {\Bbb Z}$.
Then $M_H(r,\xi+y,\chi)^{\mu}$, 
$y \in \bigoplus_{i=1}^n {\Bbb Z}C_i, \chi \in {\Bbb Z}$
is a smooth scheme of dimension
$(y^2)-2(\xi,y)+d+2r(\chi_0-\chi)+q$, $q=\dim H^1(X,{\cal O}_X)$, if
it is not empty.
\begin{lem}
\begin{enumerate}
\item
$M_H(r,\xi+x,\chi)^{\mu}$, $x \in B_{(\xi,d)}$ consists of locally free
sheaves.
\item
$H^1(C_i,E_{|C_i})=0$ for all $E \in M_H(r,\xi+x,\chi)^{\mu}$, 
$x \in B_{(\xi,d)}$.
\end{enumerate}
\end{lem}

\begin{proof}
We prove the second claim. The proof of the first one is similar.
If $H^1(C_i,E_{|C_i}) \ne 0$, then there is a surjective homomorphism
$\phi:E \to {\cal O}_{C_i}(-1-k)$, $k>0$. 
By our assumption on $E$,
$F:=\ker \phi$
is a $\mu$-stable sheaf with
$\gamma(F)=\gamma(E)-(0,C_i,-k)$. Then 
we have 
\begin{equation}
\begin{split}
M_H(\gamma(F))=&
((x-C_i)^2)-2(\xi,x-C_i)+d-2rk+q\\
<& 2r+q-2rk \leq q.
\end{split}
\end{equation}
This is impossible. Hence the claim holds.
\end{proof}

\begin{cor}\label{cor:RDP}
Let $E$ be an element of $M_H(r,\xi+x,\chi)^{\mu}$, 
$x \in B_{(\xi,d)}$. 
For a subspace $V \subset \Hom(E,{\cal O}_{C_i}(-1))$, 
$\phi:E \to V^{\vee} \otimes {\cal O}_{C_i}(-1)$ is surjective
and $\ker \phi$ is a $\mu$-stable locally free sheaf
with the Chern character $\ch(F)=\ch(E)-(\dim V)(0,C_i,0)$.
\end{cor}
We set
\begin{equation}
\begin{split}
{\frak P}_{{\cal O}_{C_i}(-1)}^{(n)}(r,\xi+x,\chi):=&
\{(E,U^{\vee})| E \in M_H(r,\xi+x,\chi)^{\mu}, U^{\vee} \subset 
\Hom(E,{\cal O}_{C_i}(-1)), \dim U=n \}\\
=& \{(E,U^{\vee})| E \in M_H(r,\xi+x,\chi)^{\mu}, U^{\vee} \subset 
\Hom(E[1],{\cal O}_{C_i}(-1)[1]), \dim U=n \}
\end{split}
\end{equation}
and define operators $e_i,f_i,h_i$. Then we have the following
which is due to Nakajima \cite[sect. 5]{Na:gauge}. 
\begin{prop}
Let ${\frak g}$ be a finite Lie algebra generated
by ${\cal O}_{C_i}(-1)$. Then 
${\frak g}$ acts on $\bigoplus_{x \in B_{(\xi,d)}}
H_*(M_H(r,\xi+x,\chi)^{\mu})$, provided the non-emptyness
of the moduli spaces.
\end{prop}

\begin{rem}
In order to compare the correspondence in 
Theorem \ref{thm:action},
we need to set ${\Bbb F}:=E[1]$ (cf. \eqref{eq:dual-system}).
\end{rem}
 
\begin{ex}
Let $\pi:X \to {\Bbb P}^1$ 
be an elliptic $K3$ surface with a section $\sigma$.
Let $f$ be a fiber of $\pi$.
Then $H:=\sigma+tf$, $t \gg 0$ is a nef and big divisor
on $X$. Let $C_i$, $i=1,2,\dots,n$ be $(-2)$-curves contracted by
$mH$. Assume that $\gcd(r,(\xi,f))=1$.
Then for $y \in \NS(X)$ with $(y,f)=0$ and  
$k \in {\Bbb Z}$,
\begin{equation}
M_H(r,\xi+y,\chi)^{\mu}
:=\left\{E \left| 
\begin{split}
& \text{ $E$ is a torsion free sheaf with $\gamma(E)=(r,\xi+y,\chi)$ }\\
& \text{ such that $E_{|\pi^{-1}(p)}$
is stable for a general point $p \in {\Bbb P}^1$}
\end{split}
\right.
\right\}.
\end{equation}
In particular, $M_H(r,\xi+y,\chi)^{\mu}$ is projective and
coincides with $M_{H'}(r,\xi+y,\chi)$
where $H'$ is an ample divisor which is
sufficiently close to $H$.
Therefore $M_H(r,\xi+y,\chi)^{\mu}$ is not empty, provided
the expected dimension is non-negative (cf. \cite{Y:5}).  
Thus all the requirements are satisfied and we have an action of
finite Lie algebra.
By Corollary \ref{cor:exist-rational}, a similar result holds for 
a rational elliptic case.
\end{ex}

\begin{rem}
Let $X \to C$ be an elliptic surface in section \ref{subsect:elliptic}.
We use the same notations.
Let ${\cal E}$ be the universal family on $Y \times X$.
Then we have a Fourier-Mukai transform
$\Phi_{X \to Y}^{{\cal E}^{\vee}}:
{\bf D}(X) \to {\bf D}(Y)$.
We set $G':=\Phi_{X \to Y}^{{\cal E}^{\vee}}({\cal O}_{\sigma})$.
Let $\xi:=\sigma+D$, $(f,D)=0$ be an effective divisor such that
$(\xi^2)=(\sigma^2)$. 
Assume that 
\begin{equation}
(\xi,x)+(x^2)<2r \text{ for all $x \in \bigoplus_{i=1}^n {\Bbb Z}C_i$}.
\end{equation}
Then $\Phi_{X \to Y}^{{\cal E}^{\vee}}$ induces an isomorphism
$M_H^G(0,\xi+x,\chi) \cong M_H^{{G'}^{\vee}}(r,\xi',\chi')$,
where $r\chi-d(\xi+x,\sigma)>0$, 
$(c_1(G'),f)=(\xi',f)$ and $((\xi+x)^2)=
({\xi'}^2)-2r\chi'-r(K_X,\xi')+r^2 \chi({\cal O}_X)$.
Moreover $\Phi_{X \to Y}^{{\cal E}^{\vee}}({\cal O}_{C_i}(-1)) 
\cong {\cal O}_{C_i'}(k_i)$
for some $k_i$ with $\chi_{{G'}^{\vee}}({\cal O}_{C_i'}(k_i))=0$. 
The action of ${\frak g}$ generated by 
$\Phi_{X \to Y}^{{\cal E}^{\vee}}({\cal O}_{C_i}(-1))$ 
is similar to the action in this section.
Indeed if $((C_i,C_j)_{i,j})$ is of type $E_8$,
then there is a divisor $D=\sum_{i=1}^8 b_i C_i'$ such that 
$r(D,C_i')=-(c_1(G'),C_i')$.
Then replacing ${\cal E}$ be
${\cal E} \otimes {\cal O}_Y(-D)$, we may assume that $k_i=-1$
for all $i>0$.   
\end{rem}

\begin{rem}
Let $Y$ be a projective surface with rational double points as singularities
and $H'$ an ample Cartier divisor on $Y$.
Assume that there is a morphism
$\phi:X \to Y$ which gives the minimal resolution and $H=\phi^{-1}(H')$.
For simplicity, we assume that there is a unique singular point
$p \in Y$. Let $Z:=\phi^{-1}(p)$ be the fundamental cycle.
For $E \in M_H(r,\xi+x,\chi)^{\mu}$,
we have an exact sequence
\begin{equation}\label{eq:ADE-cont}
0 \to E' \to E \to F \to 0
\end{equation}
such that $F$ is a successive extensions of 
${\cal O}_{C_i}(-1)$ and
$\Hom(E',{\cal O}_{C_i}(-1))=0$, $i=1,2,\dots,n$.
Then we have $E'_{|C_i} \cong 
{\cal O}_{C_i}(1)^{\oplus a_i} \oplus 
{\cal O}_{C_i}^{\oplus b_i}$.
For all $C_i$, there is an exact sequence
\begin{equation}
0 \to G \to {\cal O}_{Z} \to {\cal O}_{C_i} \to 0,
\end{equation}
where $G$ is a successive extensions of 
${\cal O}_{C_j}(-1)$ (cf. Example \ref{ex:RDP}).
Hence we see that $H^1(Z,E'_{|Z})=0$.
Then we see that $R^1 \phi_*(E')=0$.
Since $R^1 \phi_* F=0$, we get that
$\phi_* (E) \cong \phi_* (E')$ and
$R^1 \phi_*(E) \cong R^1 \phi_*(E')=0$.
Therefore we have a morphism
\begin{equation}
\begin{matrix}
\phi_* :& M_H(r,\xi+x,\chi)^{\mu}& \to & 
M_{H'}(r,\xi+x,\chi)^{\mu}\\
& E & \mapsto & \phi_*(E),
\end{matrix}
\end{equation}
where 
$M_{H'}(r,\xi+x,\chi)^{\mu}$ is the moduli space of $\mu$-stable sheaves
on $Y$.
By this morphism, we have a contraction of the Brill-Noether locus.
We can also show that $R^1 \phi_*({E'}^{\vee})=0$ and
$E'_{|Z}$ is generated by global sections.
\begin{NB}
Since $-B_{(\xi,d)}=B_{(-\xi,d)}$, we get
$R^1 \phi_*(E^{\vee})=0$ for all $E \in
M_H(r,\xi+x,\chi)^{\mu}$. In particular,
$R^1 \phi_*({E'}^{\vee})=0$ in \eqref{eq:ADE-cont}.
\end{NB}
Thus $\phi_*(E) \cong \phi_*(E')$ is a reflexive sheaf and
$E'$ is a full sheaf.
Hence the local structure of this contraction map
is an example of the studies of Ishii \cite{I:1},\cite{I:2}.
More generally,
for each moduli space $M_H(r,\xi,\chi)^{\mu}$,
let $M_H(r,\xi,\chi)^{\#}$ be the open subset consisting of $E$
such that $E$ is locally free, \eqref{eq:smooth-condition} holds
and $R^1 \phi_*(E)=R^1 \phi_*(E^{\vee})=0$.
Then we see that $H^1(C_i,E_{|C_i})=H^1(C_i,E^{\vee}_{|C_i})=0$
for all $i$ and Corollary \ref{cor:RDP} holds.
Since $R^j \phi_*({\cal O}_{C_i}(-1))=0$ for all $j$
and ${\cal E}xt^1_{{\cal O}_X}({\cal O}_{C_i}(-1),{\cal O}_X) \cong 
{\cal O}_{C_i}(-1)$,
$\ker \phi$ belongs to
$M_H(r,\xi-(\dim V) C_i,\chi)^{\#}$.
Therefore we also have similar claims for
$M_H(r,\xi,\chi)^{\#}$.
\end{rem}

\section{Equivariant sheaves}\label{sect:equiv}

In this section, we give a remark for the moduli of 
equivariant sheaves. 
Let $G$ be a finite group acting on $X$.
Let $E_0$ be an irreducible $G$-sheaf
of dimension 0, i.e.
$E_0$ does not have a non-trivial $G$-subsheaf.
Then $\Hom(E_0,E_0)^G={\Bbb C}$.

\begin{lem}\label{lem:key5}
Let $E_0$ be an irreducible $G$-sheaf
of dimension 0.
Let $E$ be a torsion free (resp. purely 1-dimensional) $G$-sheaf.
\begin{enumerate}
\item[(1)]
Then every non-trivial extension
\begin{equation}
0 \to E \to F \to E_0 \to 0
\end{equation}
defines a torsion free (resp. purely 1-dimensional) $G$-sheaf.
\item[(2)]
Let $V$ be a subspace of $\Hom(E,E_0)$.
Then  
$\phi:E \to V^{\vee} \otimes E_0$ is surjective.
Moreover,
$\ker \phi$ is a torsion free (resp. purely 1-dimensional) $G$-sheaf.
\end{enumerate}
\end{lem}
Let $H$ be a $G$-equivariant line bundle on $X$ which is ample.
\begin{defn}
A $G$-sheaf $E$ is $\mu$-stable, if
$E$ is torsion free and 
\begin{equation}
\frac{(c_1(F),H)}{\rk F}<
\frac{(c_1(E),H)}{\rk E}
\end{equation}
for all $G$-subsheaf $F$ of $E$ with $0<\rk F <\rk E$.
\end{defn}

For a $G$-sheaf $E$ on $X$, $v(E)$ denotes the class of $E$ 
in $K^G(X)$.
For a $v \in K^G(X)$, $M_H(v)^{\mu}$ is the moduli of $\mu$-stable
$G$-sheaves $E$ with $v(E)=v$.
Assume that 
\begin{equation}
\Ext^2(E,E)^G \to H^2(X,{\cal O}_X)^G
\end{equation}
is an isomorphism for all $E \in M_H(v)^{\mu}$.
We set 
\begin{equation}
\langle v(E),v(F) \rangle:=-G\text{-}\chi(E,F)=
-\sum_i (-1)^i \dim \Ext^i(E,F)^G.
\end{equation}
Let $E_1,E_2,\dots,E_s$ be a configuration of irreducible
$G$-sheaves of dimension 0
such that 
\begin{equation}
\begin{split}
 E_i \otimes K_X &\cong E_i,\\
 \Ext^1(E_i,E_i)^G &=0.\\
\end{split}
\end{equation}
Then $v(E_i)$ are $(-2)$-vectors.
We set
\begin{equation}
{\frak P}_{E_i}^{(n)}(v):=\{(E,U^{\vee})| E\in M_H(v)^{\mu},
U^{\vee} \subset \Hom(E[1],E_i[1]), \dim U=n \}
\end{equation}
and define operators $e_i,f_i,h_i$.
Then we have an action of the Lie algebra ${\frak g}$ generated by 
$v(E_1),v(E_2),\dots,v(E_s)$ on $\bigoplus_v H^*(M_H(v)^{\mu},{\Bbb C})$.

\begin{rem}
Let $X$ be an abelian surface or a K3 surface with a
symplectic $G$-action.
Assume that there is a fixed point.
By the McKay correspondence
\cite{BKR},
we have an equivalence
$\Phi:{\bf D}^G(X) \cong {\bf D}(\widetilde{X/G})$,
where $\widetilde{X/G} \to X/G$ is the minimal resolution
of $X/G$.
Hence $M_H(v)^{\mu}$ is isomorphic to a moduli space
of objects in ${\bf D}(\widetilde{X/G})$.
If $v=(r,\xi,a)$ with $\gcd(r,(\xi,H))=1$, then
$M_H(v)^{\mu}$ is projective.
Hence $M_H(v)^{\mu}$ is a holomorphic symplectic manifold
which is birationally equivalent to a moduli space $M_{H'}(w)$
of stable sheaves on $\widetilde{X/G}$,
where $w$ is the Mukai vector corresponding to $v$ via
$\Phi$.
By a result of Huybrechts \cite{Huy:1}, \cite{Huy:2}, 
there is an isomorphism $H_*( M_H(v)^{\mu},{\Bbb Z})
\cong H_*(M_{H'}(w),{\Bbb Z})$ via a convolution product
by an algebraic cycle. Hence by fixing this identification
for each $H_*(M_{H'}(w),{\Bbb Z})$,
we also have an action
of ${\frak g}$ on $\bigoplus_w H_*(M_{H'}(w),{\Bbb C})$. 
\end{rem}
 
\begin{rem}
Assume that $X={\Bbb P}^2={\Bbb C}^2 \cup \ell_{\infty}$ 
with an action of a Klein group
$G \subset SL({\Bbb C}^2)$.
Let $W$ be a $G$-vector space.
We consider the moduli of framed $G$-sheaves  
$(E,\Phi)$, where $E$ is a torsion free $G$-sheaf on
${\Bbb P}^2$ and  
$\Phi:E_{|\ell_{\infty}} \to 
{\cal O}_{\ell_{\infty}} \otimes W$ is a
$G$-isomorphism.
This is an example of Nakajima's quiver variety
and we have an action of affine Lie algebra associated to $G$
on the homology groups \cite{Na:1998}. 
In this case, we set $\langle v(E),v(F) \rangle:=
-G\text{-}\chi(E,F(-\ell_{\infty}))$ and we use the vanishing
$\Ext^2(E,E \to 
({\cal O}_{\ell_{\infty}} \otimes W \oplus E_i))=0$
to show the smoothness of ${\frak P}_{E_i}^{(n)}$.
\end{rem}

\section{Appendix}\label{sect:appendix}

\subsection{Moduli of coherent systems}\label{subsect:coh-system}
In this subsection, we shall explain how to construct the moduli space 
of coherent systems ${\frak P}_{E_i}^{(n)}(v)$. 
We start with a definition of a flat family.
\begin{defn}
Let $S$ be a scheme and
${\cal E}_{\bullet}:\cdots \to {\cal E}_{-1} \to {\cal E}_0 \to \cdots$ 
a bounded complex on $S \times X$.
\begin{enumerate}
\item
${\cal E}_{\bullet}$ is a flat family of stable complexes, if
${\cal E}_i$ are coherent sheaves on $S \times X$ which are flat over $S$
and $({\cal E}_{\bullet})_s$ are stable complexes for all
$s \in S$.
\item
$({\cal E}_{\bullet},{\cal U})$ is a family of coherent systems, if
${\cal E}_{\bullet}$ is a flat family of stable complexes
and ${\cal U}$ is a locally free subsheaf of
$\Hom_{p_S}({\cal O}_S \boxtimes E_i,{\cal E}_{\bullet})$ of rank $n$
such that ${\cal U}_s \to \Hom(E_i,({\cal E}_{\bullet})_s)$ 
is injective for all $s \in S$.
In this case, we have a resolution of $E_i$
\begin{equation}
W_{\bullet}:W_{-2} \to W_{-1} \to W_0
\end{equation}
with a morphism ${\cal U} \boxtimes W_{\bullet} \to {\cal E}_{\bullet}$
as complexes which induces the inclusion
${\cal U} \to \Hom_{p_S}({\cal O}_S \boxtimes E_i,{\cal E}_{\bullet})$.
\end{enumerate}
\end{defn}
For a quasi-isomorphism ${\cal E}_{\bullet} \to {\cal E}_{\bullet}'$
of families of stable complexes over $S$,
we take a resolution of $E_i$
\begin{equation}
W_{\bullet}:W_{-2} \to W_{-1} \to W_0
\end{equation}
such that $\Ext^p(W_j,({\cal E}_k)_s)=0$, $p>0$ for 
$j=0,-1$, $k \in {\Bbb Z}$ and all $s \in S$.
Then we see that $\Ext^p(W_{-2},({\cal E}_k)_s)=0$, $p>0$ for 
$k \in {\Bbb Z}$ and all $s \in S$.
By this choice of $W_{\bullet}$, we have an isomorphism
\begin{equation}
\Hom_{{\bf K}(S \times X)}({\cal O}_S \boxtimes W_{\bullet},
{\cal E}_{\bullet}[p]) 
\to 
\Hom_{{\bf K}(S \times X)}({\cal O}_S \boxtimes W_{\bullet},
{\cal E}_{\bullet}'[p]) 
(\cong \Ext^p({\cal O}_S \boxtimes E_i,{\cal E}_{\bullet}'))
\end{equation}
where ${\bf K}(Z)$ is the homotopy category of complexes on
$Z$.
\begin{NB}
Since $[{\cal E}_{\bullet} \to {\cal E}_{\bullet}']$ is quasi-isomorphic to 0,
\begin{equation}
\Hom_{{\bf K}(S \times X)}({\cal O}_S \boxtimes W_{\bullet},
[{\cal E}_{\bullet} \to {\cal E}_{\bullet}'][p])=0
\end{equation}
for all $p$.
\end{NB}
Hence for a family of coherent systems $({\cal E}_{\bullet}',{\cal U})$,
there is a resolution of $E_i$ 
and a family of coherent systems $({\cal E}_{\bullet},{\cal U})$ such that
we have a homotopy commutative diagram:
\begin{equation}
\begin{CD}
{\cal U} \boxtimes W_{\bullet} @>{\phi}>> {\cal E}_{\bullet}\\
@| @VVV\\
{\cal U} \boxtimes W_{\bullet} @>>> {\cal E}_{\bullet}'.
\end{CD}
\end{equation}
The choice of $\phi$ is unique, up to homotopy equivalence.
In this case, we say that $({\cal E}_{\bullet},{\cal U})$
is equivalent to $({\cal E}_{\bullet}',{\cal U})$.
 
Let $q:Q_H(v) \to M_H(v)$ be a standard $PGL(N)$-covering of
$M_H(v)$ which is an open subscheme of a suitable quot-scheme
and satisfies the following properties:
\begin{enumerate}
\item
There is a flat family of stable complexes
${\cal V}_{\bullet}:{\cal V}_{-1} \to {\cal V}_0$ on $Q_H(v) \times X$,
which is $GL(N)$-equivariant.
\item
For a flat family of stable complexes
${\cal E}_{\bullet}$ parametrized by $S$,
if we take a suitable open covering $S=\cup_{\lambda} S_{\lambda}$, then
we have a morphisms $f_{\lambda}:S_{\lambda} \to Q_H(v)$ such that
${\cal E}_{\bullet|S_{\lambda}}$ 
is quasi-isomorphic to $f_{\lambda}^*({\cal V}_{\bullet})$.
In particular  
$(q \circ f_{\lambda})_{|S_{\lambda} \cap S_{\mu}}=
(q \circ f_{\mu})_{|S_{\lambda} \cap S_{\mu}}$ and 
we have a morphism $f:S \to M_H(v)$.
\end{enumerate}
We take a locally free resolution of $E_i$ 
\begin{equation}
0 \to W_{-2} \to W_{-1} \to W_0 \to E_i \to 0
\end{equation}
such that $\Ext^p(W_j,({\cal V}_k)_t)=0$, $p>0$ for 
$j=0,-1$, $k=-1,0$ and all $t \in Q_H(v)$.
Then $\Ext^p(W_{-2},({\cal V}_k)_t)=0$, $p>0$ for 
$k=-1,0$ and all $t \in Q_H(v)$. 
We set
\begin{equation}
{\cal H}_n:=
\bigoplus_{-j+k=n}\Hom_{p_{Q_H(v)}}({\cal O}_{Q_H(v)} \boxtimes W_j,
{\cal V}_k).
\end{equation}
${\cal H}_n$, $n \in {\Bbb Z}$ are locally free sheaves on $Q_H(v)$.
We take a complex
\begin{equation}
0 \to {\cal H}_{-1} \overset{\psi_{-1}}{\to}
{\cal H}_0 \overset{\psi_0}{\to}
{\cal H}_1 \overset{\psi_1}{\to} \cdots
\end{equation}
associated to ${\bf R}\Hom_{p_{Q_H(v)}}({\cal O}_{Q_H(v)} \boxtimes E_i,
{\cal V}_{\bullet})$.
Since $\ker (\psi_{-1})_t \cong \Hom(E_i,{\cal E}_t[-1])=0$
for all $t \in Q_H(v)$,
$\psi_{-1}$ is injective as a vector bundle homomorphism.
Hence ${\cal H}_0':=\coker \psi_{-1}$ is a locally free sheaf on
$Q_H(v)$.
For the morphism $f_{\lambda}:S_{\lambda} \to Q_H(v)$ and
a locally free subsheaf ${\cal U} \subset
\Hom_{p_S}({\cal O}_S \boxtimes E_i,{\cal E}_{\bullet})$ such that
${\cal U}_s \to \Hom(E_i,({\cal E}_{\bullet})_s)$ 
is injective for all $s \in S$,
we have an inclusion as a vector bundle homomorphism:
\begin{equation}
{\cal U}_{|S_{\lambda}} \hookrightarrow 
\Hom_{p_S}({\cal O}_S \boxtimes E_i,{\cal E}_{\bullet})_{|S_{\lambda}}
=\ker (f_{\lambda}^*({\cal H}_0') \to f_{\lambda}^*({\cal H}_1))
\hookrightarrow f_{\lambda}^*({\cal H}_0'). 
\end{equation}
We take a Grassmann bundle $Gr({\cal H}_0',n) \to Q_H(v)$
over $Q_H(v)$
parametrizing $n$-dimensional subspaces $U$
of $({\cal H}_0')_t$, $t \in Q_H(v)$.
Then we have a lifting $\widetilde{f}_{\lambda}:S_{\lambda} \to
Gr({\cal H}_0',n)$
of $f_{\lambda}$ and an equivalence between 
$({\cal E}_{\bullet},{\cal U}_{|S_{\lambda}})$ and 
$(\widetilde{f}_{\lambda}^*({\cal V}_{\bullet}),
{\cal U}_{|S_{\lambda}})$.
Hence ${\frak P}_{E_i}^{(n)}(v)$ is constructed as a closed subscheme of
$Gr({\cal H}_0',n)/PGL(N)$.

\begin{NB}
Let ${\cal F}_{\bullet}$ be a flat family of complexes over over 
an affine scheme $S$. 
Assume that $\Ext^p(({\cal F}_k)_s,E_i)=0, p>0$ for all $k$ and $s \in S$.
For a complex
$W_{\bullet}':W_0' \to W_1' \to W_2'$ such that
$H^0(W_{\bullet})=E_i$ and $H^j(W_{\bullet})=0, j=1,2$,
if $\Ext^p(({\cal F}_k)_s,W_j')=0, p>0$ for $j=0,1$, $k \in {\Bbb Z}$ 
and all $s \in S$. Then
$\Ext^p(({\cal F}_k)_s,W_2')=0, p>0$ for $k \in {\Bbb Z}$ 
and all $s \in S$. Hence we have isomorphisms
\begin{equation}
\Ext^p({\cal F}_{\bullet},{\cal O}_S \otimes E_i) \cong
\Hom_{{\bf K}(S \times X)}({\cal E}_{\bullet},{\cal O}_S \otimes E_i[p])
\cong
\Hom_{{\bf K}(S \times X)}({\cal E}_{\bullet},
{\cal O}_S \otimes W_{\bullet}'[p]).
\end{equation}

\end{NB}

\subsection{The existence of stable sheaves on a
rational elliptic surface}\label{subsect:exist}


We shall find the conditions for the existence of stable
sheaves on a rational elliptic surface 
$\pi:X \to {\Bbb P}^1$ 
with a section $\sigma$. 
We first note that a divisor $C$ with $(C^2)=(C,K_X)=-1$
is effective.
Indeed since $(K_X-C,f)=-1$, $H^2(X,{\cal O}_X(C))=0$.
By the Riemann-Roch theorem,
$\dim H^0(X,{\cal O}_X(C)) \geq \chi({\cal O}_X(C))=1$.
The following is the result for the case of rank 0.

\begin{prop}\label{prop:rational-exist}
Let $X$ be a rational elliptic surface with a section $\sigma$. 
Let $D$ be a divisor with $(D^2) \geq 0$.
Assume that $(0,D,\chi)$ is primitive.
Then $M_H^G(0,D,\chi)$ is not empty for a general $H$ and $G$
if and only if
$(D,C) \geq 0$ for all divisor $C$ with
$(C^2)=(C,K_X)=-1$.
\end{prop}

\begin{proof}
We use the notation in subsection \ref{subsect:rational}.
Since $\overline{M}_{({\cal X},{\cal H})/T}^G(0,D,\chi) \to T$ is smooth,
it is sufficient to prove the claim for a nodal rational elliptic
surface $X$.
Let $C$ be a divisor with $(C^2)=(C,K_X)=-1$.
Since every fiber is irreducible, $C$ must be a section of $\pi$. 
If $(D,C)<0$, then
$\chi({\cal O}_{C}(k),E)=-(D,C)>0$ for
all sheaf $E$ with $c_1(E)=D$.
We set $n:=\max \{k|\Hom({\cal O}_{C}(k),E) \ne 0\}$.
Then $\Hom({\cal O}_{C}(n),E) \ne 0$ and
$\Hom(E,{\cal O}_{C}(n))^{\vee}=\Ext^2({\cal O}_{C}(n+1),E) \ne 0$.
This means that $E$ is not semi-stable, unless 
$E \cong {\cal O}_C(n)$.

Conversely, we assume that $(D,C) \geq 0$ for all
section $C$ with $(C^2)=(C,K_X)=-1$. 
Then $D$ is a nef divisor.  
If $(D,f)=1$, then there is a section $\tau$ of $\pi$ such that
$D=\tau+nf$, $n > 0$.
In this case, $M_H(0,\tau+nf,\chi) \cong \Hilb_X^n \ne \emptyset$
via the relative Fourier-Mukai transform.
Since the non-emptyness does not depend on
the choice of $G$ \cite{Y:11},
we get our claim.
Hence we may assume that $(D,f) \geq 2$.
We shall show that there is a reduced and irreducible curve
$C \in |D|$. Then a line bundle $E$ on $C$ with
$\chi(E)=\chi$ belongs to 
$M_H(0,D,\chi)$. 
If $(D^2) \geq 1$ or $(D,f) \geq 3$, 
then $D':=D-K_X$ is a nef divisor with
$({D'}^2) \geq 5$.
Assume that there is an effective divisor $B$ with
$(D',B) \leq 1$. Since $0 \leq (D,B) \leq (D',B) \leq 1$,
(i) $(f,B)=0$ and $(D,B)\leq 1$ or (ii) $(f,B)=1$ and $(D,B)=0$.
In the first case, $B=nf$.
Since $(D,f) \geq 2$, this is impossible.
In the second case, there is a section $\tau$ and $B=\tau+nf$. 
Then $(B^2)=2n-1 \ne 0$.
By the Reider's result \cite{R}, $D=D'+K_X$ is base point free.

If $(D^2)=0$ and $(D,f)=2$, 
then $D=2\tau_1+f$ or $D=\tau_1+\tau_2$ with
$(\tau_1,\tau_2)=1$, where $\tau_1,\tau_2$ are sections of $\pi$.
In the first case, $(D,\tau_1)=-1$, which is a contradiction.
In the second case, $D$ is connected and 
$D$ is base point free.
By Bertini's theorem, there is a reduced and irreducible curve
$C \in |D|$.
\end{proof}

\begin{defn}
We set 
\begin{equation}
{\cal C}:=\left\{D \in \Pic(X)\left|
\begin{split}
& \;(D,C) \geq 0
\text{ for all divisors $C$ }\\
& \text{ with $(C^2)=(C,K_X)=-1$}
\end{split}
\right.\right\}.
\end{equation}
Let $W:=W(E_8^{(1)})$ be the Weyl group
of the sublattice $f^{\perp} \cong E_8^{(1)}$ of $\Pic(X)$.
$W$ acts on $\Pic(X)$ and
${\cal C}$ is a $W$-invariant subset of $\Pic(X)$. 
Let ${\cal C}^+ \subset {\cal C}$ be the set of nef divisors. 
If $X$ is nodal, then ${\cal C}^+={\cal C}$.
 \end{defn}

\begin{thm}\label{thm:rational-exist}
Let $r$ and $d$ be relatively prime integers with $r \geq 0$.
\begin{enumerate}
\item
For any $D \in \langle \sigma,f \rangle^{\perp}$,
there is a stable vector bundle $E_D$ such that
$\rk(E_D)=r$, $c_1(E_D) \equiv d\sigma+D \mod {\Bbb Z}f$ and
$\chi(E_D,E_D)=1$. $E_D$ is unique up to $E_D(nf)$, $n \in {\Bbb Z}$.
We set 
\begin{equation}
{\cal E}(r,d):=\{E_D | (D,\sigma)=(D,f)=0 \}.
\end{equation}
\item
Let $F \in K(X)$ be a primitive class with
$\rk(F)=lr$ and $(c_1(F),f)=ld$.
Assume that $\chi(F,F) \leq 0$.
We take an ample divisor $H$ which is sufficiently close to $f$.
Then $F$ is represented by a stable sheaf if and only if
$\chi(E_D,F) \leq 0$ for all $E_D \in {\cal E}(r,d)$.
Moreover $F$ is represented by a $\mu$-stable vector bundle,
if $lr>1$.
\end{enumerate}
\end{thm}

\begin{proof}
We may assume that $lr>0$.
By the deformation argument in the proof of
Proposition \ref{prop:rational-exist}, 
we may assume that $X$ is nodal.
We first prove (i).
We note that $M_H(0,rf,-d) \cong X$.
Let ${\cal E}$ be a universal family
on $X \times X$. 
Since every fiber is irreducible,
we have $\sigma-D=\tau-((\sigma,\tau)+1)f$, 
where $\tau$ is a section of $\pi$.
Then ${\cal E}^{\vee}_{|X \times \tau}$ is a stable sheaf
with the desired invariant.
We next prove (ii).
The proof of the necessary condition is similar to 
the proof of Proposition \ref{prop:rational-exist}.
We shall show that the condition is sufficient. 
Let $\Phi_{X \to X}^{\cal E}:{\bf D}(X) \to {\bf D}(X)$ 
be the relative Fourier-Mukai transform defined by the sheaf ${\cal E}$.
Then $\Phi_{X \to X}^{\cal E}(E_D)[1]={\cal O}_{\tau}$, where
$\tau$ is a section of $\pi$ such that $\tau-\sigma 
\equiv -D \mod {\Bbb Z}f$.
Then $\rk(\Phi_{X \to X}^{\cal E}(F)[1])=0$ and 
$c_1(\Phi_{X \to X}^{\cal E}(F)[1]) \in {\cal C}$.
Therefore $\Phi_{X \to X}^{\cal E}(F)[1]$ 
is represented by a line bundle $L$ on
a reduced and irreducible curve.
Then the inverse $\Phi_{X \to X}^{{\cal E}^{\vee}}(L)[1]$
is a $\mu$-stable sheaf.
\end{proof}

By the proof of the theorem, we also get the following.

\begin{cor}\label{cor:exist-rational}
If $\gcd(r,(\xi,f))=1$ and the expected dimension is non-negative, 
then $M_H(r,\xi,\chi)$ is not empty,
where $H$ is sufficiently close to $f$.
\end{cor}

Let $X$ be a rational elliptic surface with a section $\sigma$
such that there is a singular fiber $\pi^{-1}(o)=\sum_{i=0}^8 a_i C_i$,
$o \in {\Bbb P}^1$ of type $E_8^{(1)}$,
where $C_i$ are smooth $(-2)$-curves.
We assume that $a_0=1$.
Let $C$ be a divisor with $(C^2)=(C,K_X)=-1$.
Then $C=\sigma+\sum_{i=0}^8 n_i C_i$, $n_i \geq 0$. 
Hence
\begin{equation}
{\cal C}^+=\left\{D \in \Pic(X) \left|
 (D,\sigma) \geq 0, (D,C_i) \geq 0, 0 \leq i \leq 8 \right. \right\}.
\end{equation}
Thus $D:=r\sigma+nf+\xi$, $\xi \in \bigoplus_{i=1}^8 {\Bbb Z}C_i$.
is nef if and only if 
\begin{equation}
\begin{cases}
n \geq r,\\
(\xi,C_i) \geq 0,\;1 \leq i \leq 8\\
\sum_{i=1}^8 a_i (\xi,C_i) \leq r.
\end{cases}
\end{equation}
Let $W$ be the affine Weyl group of $E_8^{(1)}$.
Then $M_H(0,D',\chi) \ne \emptyset$
if and only if $D'=w(D)$ with $D \in {\cal C}^+, w \in W$.

\begin{NB}
As in \cite{GKV},
we have similar action of the quantum groups
over the function spaces of 
the set of ${\Bbb F}_q$-rational points of $M_H(v)$. 
We assume that $X$ and $M_H(v)$ are defined over a finite field
${\Bbb F}_q$ of $q$-elements.
Let $\Fun(M_H(v))$ be the set of functions over
the set of ${\Bbb F}_q$-rational points of $M_H(v)$.
We assume that all the required condions for $M_H(v)$ and $E_i$
are satisfied. 
For a pair $(a,b)$ with $b-a=\langle v_i,v \rangle$,
we set
\begin{equation}
\frak{P}_{E_i}^{(1)}(v,a,b):=
\{({\Bbb E},U) \in \frak{P}_{E_i}^{(1)}(v)| \dim \Hom(E_i,{\Bbb E})=a,
\dim \Ext^1(E_i,{\Bbb E})=b \}.
\end{equation} 
Let 
$C_{E_i}(v,a,b)$ be the characteristic function of 
$\frak{P}_{E_i}^{(1)}(v,a,b)$.
We define a function $K_{E_i}$ 
on $M_H(v) \times M_H(v-v_i)$
supported on $\frak{P}_{E_i}^{(1)}(v)$ as
\begin{equation}
\begin{split}
K_{E_i}:=&\sum_b q^{-\frac{b}{2}}C_{E_i}(v,a,b).
\end{split}
\end{equation}
We fefine operators by the convolution product
of functions:
\begin{equation}
\begin{matrix}
e_i:&\Fun(M_H(v)) & \to &  \Fun(M_H(v-v_i))\\
& f & \mapsto & p_{1*}(K_{E_i}\cdot p_2^*(f)),\\
\end{matrix}
\end{equation}

\begin{equation}
\begin{matrix}
f_i:&\Fun(M_H(v-v_i)) & \to &  \Fun(M_H(v))\\
& f & \mapsto & p_{2*}(q^{\frac{\langle v_1,v \rangle+1}{2}}
K_{E_i}\cdot p_1^*(f)),
\end{matrix}
\end{equation}

\begin{equation}
\begin{matrix}
k_i:&\Fun(M_H(v)) & \to &  \Fun(M_H(v))\\
& f & \mapsto & q^{\frac{\langle v_1,v \rangle}{2}}f.
\end{matrix}
\end{equation}

Since the convolution product satisfies
\begin{equation}
\begin{split}
(\omega(C_{E_i}(v,a,b))*C_{E_i}(v,a',b'))({\Bbb E},{\Bbb E'})
=&
\begin{cases}
\frac{q^a-1}{q-1}, & a=a',\;{\Bbb E}={\Bbb E}'\\
1, & a=a',\;{\Bbb E} \ne {\Bbb E}'\\
0, & a \ne a'
\end{cases},\\
(C_{E_i}(v+v_i,a,b)*\omega(C_{E_i}(v+v_i,a',b')))({\Bbb E},{\Bbb E'})
=&
\begin{cases}
\frac{q^b-1}{q-1}, & a=a',\;{\Bbb E}={\Bbb E}'\\
1,& a=a',\;{\Bbb E}\ne {\Bbb E}'\\
0, & a \ne a'
\end{cases},
\end{split}
\end{equation}
we see that
\begin{equation}
\begin{split}
[e_i,f_j]=&
\delta_{i,j}\frac{k_i-k_i^{-1}}{q^{\frac{1}{2}}-q^{-\frac{1}{2}}},\\
k_i f_j k_i^{-1}=& q^{\frac{\langle v_i,v_j \rangle}{2}}f_j,\\
k_i e_j k_i^{-1}=& q^{\frac{-\langle v_i,v_j \rangle}{2}}f_j.
\end{split}
\end{equation}
By the integrability of the representation,
we also have the quantum Serre relations.
Thus we have an action of $U_{q^{\frac{1}{2}}}({\frak g})$.

\end{NB}

{\it Acknowledgement.}
I would like to thank Hiraku Nakajima for valuable discussions on this subject
for years.

\end{document}